\documentclass{amsart}
\usepackage{amsthm}
\usepackage{amsmath}
\usepackage{amsfonts}
\usepackage{amssymb}
\usepackage{graphicx}
\usepackage{float}
\usepackage{rotating}
\restylefloat{figure}
\usepackage{mathrsfs}
\usepackage{fancyhdr}
\usepackage[breaklinks]{hyperref}
\usepackage[numbers,sort&compress]{natbib}
\theoremstyle{plain}
\newtheorem{lem}{Lemma}
\newtheorem{prop}[lem]{Proposition}
\newtheorem{thm}[lem]{Theorem}

\theoremstyle{definition}

\newtheorem{rem}[lem]{Remark}

\newcommand{\R}{\mathbb R}

\newcommand{\N}{\mathbb N}
\newcommand{\dx}{\,dx}
\newcommand{\dz}{\,dz}

\newcommand{\eps}{\varepsilon}
\expandafter\mathchardef\expandafter\varphi\number\expandafter\phi\expandafter\relax
\expandafter\mathchardef\expandafter\phi\number\varphi
\newcommand{\norm}[1]{\left\Vert#1\right\Vert}

\renewcommand{\autoref}[1]{\text{Eq.}~\eqref{#1}}
\newcommand{\bea}{\begin{eqnarray}}
\newcommand{\eea}{\end{eqnarray}}
\newcommand{\beq}{\begin{equation}}
\newcommand{\eeq}{\end{equation}}
\begin{document}
\title{On an elliptic-parabolic MEMS model with two free boundaries}
\author{Martin Kohlmann}
\address{Dr.\ Martin Kohlmann, Goerdelerstraße 36, 38228 Salzgitter, Germany}
\email{martin\_kohlmann@web.de}
\keywords{MEMS, free boundary problem, local and global well-posedness, non-existence, asymptotic stability, small aspect ratio limit}
\subjclass[2010]{35R35, 35M33, 35B30, 35Q74, 74M05}
\date{\today}
\begin{abstract} We discuss an evolution free boundary problem of mixed type with two free boundaries modeling an idealized electrostatically actuated MEMS device. While the electric potential is the solution of an elliptic equation, the dynamics of the membranes' displacement is modeled by two parabolic equations. It is shown that the model is locally well-posed in time and that solutions exist globally for small source voltages whereas non-existence holds for large voltage values. Moreover, our model possesses a steady state solution that is asymptotically stable. Finally, we show that in the vanishing aspect ratio limit, solutions of the model converge towards solutions of the associated small aspect ratio problem.
\end{abstract}
\maketitle
\section{Introduction and main results}\label{sec_intro}
\begin{figure}[H]
\begin{center}
\includegraphics[width=8cm]{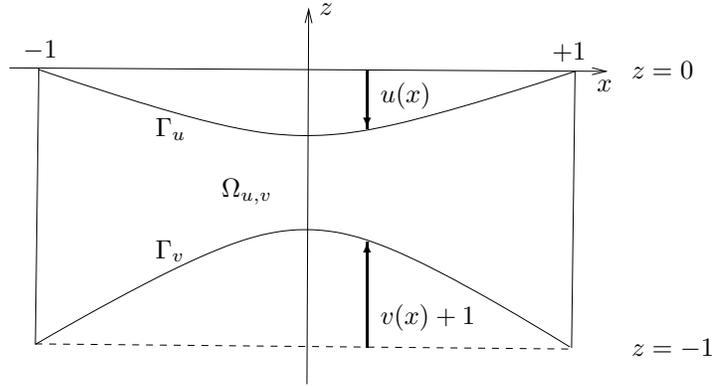}
\end{center}
\caption{An idealized model for an electrostatic MEMS device with two free boundaries.}
\label{figabstr}
\end{figure}

Mathematical models for Micro-Electro Mechanical Systems (MEMS) have been studied with regularity in the last few years, cf.~\cite{PB03} for an overview and \cite{GG0607,GG08,EGG10,WL12} for some more recent references. In \cite{mk13} the stationary version of the following free boundary problem for an idealized electrostatic MEMS device has been proposed: Let $I=(-1,1)$, pick $\tau>0$ and $q\in(2,\infty)$, consider functions $u,v\in C([0,\tau);\,W_q^2(I))\cap C^1([0,\tau);\,L_q(I))$ with $-1\leq v<u\leq 0$ on $I$, let
$$\Omega_{u(t),v(t)}=\{(x,z)\in I\times(0,-1);\,v(t,x)<z<u(t,x)\}$$
and denote by $\Gamma_{u(t)}=\{z=u(t,x)\}$ and $\Gamma_{v(t)}=\{z=v(t,x)\}$ the horizontal boundary components of $\Omega_{u(t),v(t)}$. The functions $u$ and $v$ model the one-dimensional displacements of two deformable elastic membranes from $\Gamma_{0}$ and $\Gamma_{-1}$ when a non-zero source voltage is applied to the device $\Omega_{0,-1}$; see Fig.~\ref{figabstr}. Since both membranes should be held fixed along the boundary of the device, we impose the conditions $u(t,\pm 1)=0$ and $v(t,\pm 1)=-1$. The evolution of the membranes starts from $u(0,x)=u_0(x)$ and $v(0,x)=v_0(x)$. The electrostatic potential $\phi$ in the region between both membranes satisfies the Laplace equation, is equal to zero on the lower and one on the upper membrane and is a linear function of $z$ on $\{x=\pm 1\}\cap\overline\Omega_{u(t),v(t)}$. Moreover, the functions $u$ and $v$ solve a heat equation with a right-hand side proportional to the square of the trace of the gradient of the electrostatic potential on the respective membrane. From the modeling point of view, we also need two parameters $\lambda,\mu>0$ in the equations on the free boundaries proportional to the square of the source voltage and inversely proportional to the surface tension of the respective membrane. The coefficients $\lambda$ and $\mu$ interrelate the strengths of the electrostatic and mechanical forces in the device. Finally, by nondimensionalization, there is a parameter $\eps>0$ called the \emph{aspect ratio} of the device, comparing gap size to device length.
\begin{picture}(0,0)
\put(135,128){$\Omega_{u,v}$}
\put(110,105){$\Gamma_v$}
\put(110,151){$\Gamma_u$}
\put(190,176){\vector(0,-1){22}}
\put(190,71){\vector(0,1){40}}
\put(195,164){$u(x)$}
\put(195,80){$v(x)+1$}
\put(290,173){$z=0$}
\put(290,68){$z=-1$}
\put(260,180){$+1$}
\put(60,181){$-1$}
\put(172,197){$z$}
\put(277,168){$x$}
\end{picture}
Let $\partial_t=\frac{\partial}{\partial t}$, $\partial_x=\frac{\partial}{\partial x}$, $\partial_z=\frac{\partial}{\partial z}$, $\nabla_\eps=(\eps\partial_x,\partial_z)$ and $\Delta_\eps=\eps^2\partial_x^2+\partial_z^2$. Our problem reads
\begin{align}
-\Delta_\eps\phi & = 0, & \text{in }\Omega_{u,v},\,t>0,\label{originalproblem1}\\
\phi & = \frac{z-v}{u-v}, & \text{on }\partial\Omega_{u,v},\,t>0,\label{originalproblem2}\\
\partial_tu-\partial_x^2u & = -\lambda|\nabla_\eps\phi|^2,&\text{on }\Gamma_u, \,t>0,\label{originalproblem3}\\
\partial_tv-\partial_x^2v & = \mu|\nabla_\eps\phi|^2,&\text{on }\Gamma_v,\,t>0,\label{originalproblem4}\\
u(t,\pm 1) & = 0,&t>0,\label{originalproblem5}\\
v(t,\pm 1) & =-1,&t>0,\label{originalproblem6}\\
u(0,x) & = u_0,&x\in I,\label{originalproblem7}\\
v(0,x) & = v_0,&x\in I.\label{originalproblem8}
\end{align}
In physics, the problem \eqref{originalproblem1}--\eqref{originalproblem8} serves as a model for a so-called DFM device, i.e., a MEMS with double freestanding membranes as explained in, e.g., \cite{FTYS13}. For the convenience of the reader, a derivation of the model \eqref{originalproblem1}--\eqref{originalproblem8} can be found in the Appendix. For time-independent functions $(u,v,\phi)$, the system \eqref{originalproblem1}--\eqref{originalproblem8} reduces to the problem studied in \cite{mk13} where the existence of solutions in suitable Sobolev and H\"older spaces for small source voltages (i.e., small values of the parameters $\lambda,\mu$) has been proved. The results of the paper at hand refer to two of the open problems stated in \cite{mk13} which is why this work can be seen as a companion paper of \cite{mk13}.

Experience has shown that studying the \emph{small aspect ratio limit} ($\eps\to 0$) of an idealized MEMS model is useful for getting results on the existence and uniqueness of solutions \cite{FMP03,FMPS0607,EGG10,GG0607,GG08,H11,GHW09,G08a,G08b,LY07,PB03}. Sending $\eps\to 0$, one obtains the following narrow gap model from \eqref{originalproblem1}--\eqref{originalproblem8}:
\begin{align}
\phi & = \frac{z-v}{u-v}, & \text{in }\overline{\Omega}_{u,v},\,t>0,\label{sarm1}\\
\partial_tu-\partial_x^2u & = -\frac{\lambda}{(u-v)^2},&x\in I,\, t>0,\label{sarm2}\\
\partial_tv-\partial_x^2v & = \frac{\mu}{(u-v)^2},&x\in I,\, t>0,\label{sarm3}\\
u(t,\pm 1) & = 0,&t>0,\label{sarm4}\\
v(t,\pm 1) & =-1,&t>0,\label{sarm5}\\
u(0,x) & = u_0,&x\in I,\label{sarm6}\\
v(0,x) & = v_0,&x\in I.\label{sarm7}
\end{align}
The right-hand side of \eqref{sarm2} and \eqref{sarm3} has a singularity for $u(x)=v(x)$. This singularity corresponds to the physical observation that both membranes come closer and closer and finally touch when the source voltage is increased. This phenomenon called \emph{pull-in instability} is a major factor limiting the effectiveness of many real-life MEMS devices. Thus for practical reasons it is important to know the precise value of the pull-in voltage such that there is a stable configuration of the device below the threshold and collision of the membranes and malfunction for voltages larger than or equal to the threshold.
Since the parameters $\lambda$ and $\mu$ are proportional to the source voltage of the idealized MEMS device $\Omega_{u(t),v(t)}$, it is reasonable to expect that our model possesses solutions for small values of $\lambda$ and $\mu$ which cease to exist as $(\lambda,\mu)$ vary though the parameter space.

There are a wide range of papers to comment on where the authors suppose that the idealized MEMS device consists of only one free membrane that is suspended above a rigid, fixed ground plate and is held fixed along the boundary. We refer the reader to \cite{PT00,P01,P0102,PB03} for a detailed presentation of some important results for this type of model. In \cite{WL12,ELW12} the authors study the stationary and the dynamic free boundary problem associated with the model with a fixed ground plate. In \cite{ELW13} an elliptic-parabolic problem with an additional curvature term is discussed. Recently, some fourth-order models including the mechanical effects damping and bending have been studied. In this case, terms of the form $\alpha\partial_t^2u$ and $\beta\partial_x^4u$ occur in the equation on the free boundary, cf.~also \cite{WL13,LL12}. Hyperbolic ($\alpha>0$) MEMS models associated with a device with only one free membrane have been subject of \cite{FMP03,G10,KLNT11} and further references concerning second-order parabolic ($\alpha=\beta=0$) models are, e.g., \cite{FMPS0607,GG08,G08a,G08b,H11}.

Following a line of arguments of Lauren\c{c}ot's paper \cite{ELW12}, the results of the present paper and its organization are as follows: In Section~\ref{sec_lwp} we show that the problem \eqref{originalproblem1}--\eqref{originalproblem8} is locally well-posed for any pair of values $\lambda,\mu>0$. To this end, the free boundary problem \eqref{originalproblem1}--\eqref{originalproblem8} is mapped to a reference problem of mixed type on a fixed domain. Solving the elliptic equation for the potential first, our analysis results in a semilinear evolution equation for the free surfaces with a right-hand side depending on the trace of the gradient of the potential. We then apply the Contraction Mapping Theorem to obtain a solution $(u_\eps,v_\eps,\phi_\eps)$ of regularity $W_q^2(I)^2\times W_2^2(\Omega_{u,v})$, for any $\eps>0$. Furthermore, it is shown that this solution exists globally in time if $\lambda,\mu<m_{1}$, for some $m_1(\eps)>0$. Section~\ref{sec_nonex} deals with non-existence of global solutions. Using a suitable Lyapunov functional, we compute a number $m_2(\eps)>0$ such that, for $\max\{\lambda,\mu\}>m_2$, the maximal existence time of the solution to \eqref{originalproblem1}--\eqref{originalproblem8} is finite. A smooth branch of steady state solutions of \eqref{originalproblem1}--\eqref{originalproblem8} emanating from $(\lambda,\mu)=(0,0)$ is obtained in Section~\ref{sec_stab} from the Implicit Function Theorem. Applying the Principle of Linearized Stability, we also show that this steady state is asymptotically stable. Finally, in Section~\ref{sec_sar}, a rigorous justification of the small aspect ratio model \eqref{sarm1}--\eqref{sarm7} is given by showing that there is $\eps_*>0$ so that any family of solutions $\{(u_\eps,v_\eps,\phi_\eps);\,\eps<\eps_*\}$ to \eqref{originalproblem1}--\eqref{originalproblem8} contains a sequence that converges to a solution of \eqref{sarm1}--\eqref{sarm7} in suitable spaces. A discussion of our results can be found in Section~\ref{sec_disc} and the Appendix contains the derivation of our model from the physical viewpoint.
\section{Local and global well-posedness}\label{sec_lwp}
As major difference to the stationary version of \eqref{originalproblem1}--\eqref{originalproblem8} studied in \cite{mk13}, we will use a $W^2_q$-setting, $q\in(2,\infty)$, to be able to work with the heat semigroup in $L_q(I)$. 
We introduce, for $q\in(2,\infty)$ and $\kappa\in(0,1/2)$, the sets
\begin{align}
S_q(\kappa)&:=\bigg\{(u,v)\in W^2_{q}(I)\times W^2_{q}(I);\,(u,v)(\pm 1)=(0,-1),\,\norm{u}_{W^2_{q,D}(I)}<\frac{1}{\kappa},\nonumber\\
& \hspace{1cm}\norm{v+1}_{W^2_{q,D}(I)}<\frac{1}{\kappa},\, v(x)+\kappa<u(x)-\kappa,\,\forall x\in I\bigg\}\nonumber
\end{align}
where
$$
W^{2\alpha}_{q,D}(I):=\left\{
\begin{array}{ll}
\{w\in W^{2\alpha}_{q}(I);\,w(\pm 1)=0\}, & 2\alpha\in(1/q,2], \\
W_q^{2\alpha}(I), & 0\leq 2\alpha<1/q.
\end{array}
\right.
$$
%
Then $S_q(\kappa)+\{(0,1)\}$ is open in the topology of $W^2_{q,D}(I)\times W^2_{q,D}(I)$ and its closure $\overline S_q(\kappa)$ is obtained by replacing $<$ by $\leq$ in the definition of $S_q(\kappa)$.

In a first step, we transform the problem \eqref{originalproblem1}--\eqref{originalproblem8} on the a priori unknown domain $\Omega_{u(t),v(t)}$ to the fixed reference domain $\Omega:=I\times(0,1)$ by applying a time-dependent transformation of coordinates $T=T_{u(t),v(t)}\colon\overline{\Omega}_{u(t),v(t)}\to\overline{\Omega}$ given by
\beq\label{defT}T(x,z)=(x',z')=\left(x,\frac{z-v(t,x)}{u(t,x)-v(t,x)}\right).\eeq
It is easily checked that, with the definition of $\overline{\Omega}_{u(t),v(t)}$ in Section~\ref{sec_intro}, $T_{u(t),v(t)}$ is a diffeomorphism $\overline{\Omega}_{u(t),v(t)}\to\overline{\Omega}$ with the inverse
$$T^{-1}(x',z')=(x',z'(u(t,x')-v(t,x'))+v(t,x')).$$
Let $\theta^*(u,v)$ and $\theta_*(u,v)$ be the pull-back and push-forward operators for the pair $(\Omega_{u,v},\Omega)$ defined by $\theta^*(u,v)\tilde w=\tilde w\circ T_{u,v}$ and $\theta_*(u,v)w=w\circ T^{-1}_{u,v}$ where $w$ and $\tilde w$ are functions of the coordinates $(x,z)$ and $(x',z')$ respectively, i.e.,
$$[\theta^*(u,v)\tilde w](x,z)=\tilde w(T_{u,v}(x,z))\quad\text{and}\quad [\theta_*(u,v) w](x',z')=w(T_{u,v}^{-1}(x',z')).$$
We let $\widetilde{\Delta}_{u,v;\eps}=\theta_*(u,v)\Delta_\eps\theta^*(u,v)$ denote the time-dependent transformed Laplace operator on $\Omega$. As explained in \cite{mk13},
\begin{align}\widetilde{\Delta}_{u,v;\eps}\tilde w & = \eps^2\tilde w_{x'x'} - 2\eps^2\tilde w_{x'z'}\frac{z'(u_{x'}-v_{x'})+v_{x'}}{u-v}
+\tilde w_{z'z'}\frac{1+\eps^2[z'(u_{x'}-v_{x'})+v_{x'}]^2}{(u-v)^2} \nonumber\\
& \quad +\eps^2\tilde w_{z'}\left(2\frac{u_{x'}-v_{x'}}{(u-v)^2}[z'(u_{x'}-v_{x'})+v_{x'}]-\frac{z'(u_{x'x'}-v_{x'x'})+v_{x'x'}}{u-v}\right);\label{deltatilde}
\end{align}
here the notation $u_{x'}$ stands for $\partial_{x'}u$ et cetera. We first concentrate on the elliptic boundary value problem \eqref{originalproblem1}--\eqref{originalproblem2} which is reformulated as
\begin{align}
-\left(\widetilde{\Delta}_{u(t),v(t);\eps}\tilde\phi\right)(t,x',z') & = 0, & \hspace{-2cm} (x',z')\in\Omega,\, t>0, \label{phieqn1}\\
\tilde\phi (t,x',z') & = z', & \hspace{-2cm} (x',z')\in\partial\Omega,\, t>0, \label{phieqn2}
\end{align}
with $\tilde\phi=\theta_*(u(t),v(t))\phi$. With $\psi(t,x',z')=\tilde\phi(t,x',z')-z'$ and
\beq f_{u,v;\eps} = \eps^2\left(2\frac{u_{x'}-v_{x'}}{(u-v)^2}[z'(u_{x'}-v_{x'})+v_{x'}]-\frac{z'(u_{x'x'}-v_{x'x'})+v_{x'x'}}{u-v}\right)\label{deff}\eeq
we can rewrite the problem \eqref{phieqn1}--\eqref{phieqn2} as
\begin{align}
-\left(\widetilde{\Delta}_{u(t),v(t);\eps}\psi\right)(t,x',z') & = f_{u(t),v(t);\eps}, & \hspace{-1cm} (x',z')\in\Omega,\, t>0, \label{phieqn1'}\\
\psi (t,x',z') & = 0, & \hspace{-1cm} (x',z')\in\partial\Omega,\, t>0. \label{phieqn2'}
\end{align}
In the following, $c_1,c_2,c_3,\ldots$ stand for positive constants depending on what is postpositioned in brackets. For $(u,v)\in{S}_q(\kappa)$ and $x\in I$, we have
$$u(x)-v(x)>2\kappa,\quad \norm{u}_{C^1([-1,1])},\norm{v}_{C^1([-1,1])}\leq c_1(\kappa,q).$$
Then it easy to see that $-\widetilde{\Delta}_{u(t),v(t);\eps}$ is strictly elliptic, with an ellipticity constant independent of $(u,v)\in S_q(\kappa)$, and writing \eqref{deltatilde} in divergence form,
\begin{align}
-\widetilde{\Delta}_{u,v;\eps}\tilde w & = -\partial_{x'} \left(\eps^2\tilde w_{x'}-\eps^2\frac{z'(u_{x'}-v_{x'})+v_{x'}}{u-v}\tilde w_{z'}\right) \nonumber\\
& \quad - \partial_{z'} \left(-\eps^2\frac{z'(u_{x'}-v_{x'})+v_{x'}}{u-v}\tilde w_{x'} + \frac{1+\eps^2[z'(u_{x'}-v_{x'})+v_{x'}]^2}{(u-v)^2}\tilde w_{z'}\right) \nonumber\\
& \quad +\eps^2(u_{x'}-v_{x'})\frac{z'(u_{x'}-v_{x'})+v_{x'}}{(u-v)^2}\tilde w_{z'}-\eps^2\frac{u_{x'}-v_{x'}}{u-v}\tilde w_{x'},\nonumber\\
& = -\partial_{x'}(a_{11}(u,v;\eps)\tilde w_{x'}+a_{12}(u,v;\eps)\tilde w_{z'}) + b_1(u,v;\eps)\tilde w_{x'}\nonumber\\
&\quad -\partial_{z'}(a_{21}(u,v;\eps)\tilde w_{x'}+a_{22}(u,v;\eps)\tilde w_{z'}) + b_2(u,v;\eps)\tilde w_{z'},\nonumber
\end{align}
it is clear that
$$\sum_{i,j=1}^2\norm{a_{ij}(u,v;\eps)}_{W^1_q(\Omega)}+\sum_{i=1}^2\norm{b_i(u,v;\eps)}_{L_\infty(\Omega)}\leq c_2(\kappa,\eps),\quad\forall(u,v)\in S_q(\kappa),$$
and that $a_{ij}$, $1\leq i,j\leq 2$, belongs to $C(\overline\Omega)$. Since $f_{u,v;\eps}\in L_2(\Omega)$, we can apply the arguments in the proofs of Proposition~2.1 and Lemma~2.2 of \cite{ELW12} to obtain that the problem \eqref{phieqn1'}--\eqref{phieqn2'} possesses, for $q\in(2,\infty)$, $\kappa\in(0,1/2)$,  $\eps>0$ and $(u,v)\in S_q(\kappa)$, a unique solution $\psi_{u,v;\eps}\in W^2_{2,D}(\Omega)$  satisfying
\beq\label{bddinv}\norm{\psi_{u,v;\eps}}_{W^2_{2,D}(\Omega)}\leq c_3(\kappa,\eps)\norm{f_{u,v;\eps}}_{L_2(\Omega)}.\eeq
Hence $\tilde\phi_{u,v;\eps}=\psi_{u,v;\eps}+z'$ is the unique solution to \eqref{phieqn1}--\eqref{phieqn2} with $\norm{\tilde\phi_{u,v;\eps}}_{W^2_2(\Omega)}\leq c_4(\kappa,\eps)$. Moreover, with the notation $\tilde w(x'):=w(-x')$, $x'\in I$, we have that $\tilde\phi_{\tilde u,\tilde v;\eps}(t,x',z')=\tilde\phi_{u,v;\eps}(t,-x',z')$ for $(x',z')\in\Omega$, as $S_q(\kappa)$ is invariant under the operation $x\mapsto-x$ and by \eqref{deltatilde} and uniqueness. In particular, for even functions $(u,v)$, the potential $\tilde\phi$ is even in $x'$.

Next, we discuss the parabolic equations \eqref{originalproblem3}--\eqref{originalproblem4} on the free boundaries. We first concentrate on the right-hand side terms of these equations. For simplicity, we write $x$ instead of $x'$ henceforth. For $(u,v)\in S_q(\kappa)$, let
\begin{align}
g_\eps(u,v)=\left(\frac{1+\eps^2u_{x}^2}{(u-v)^2}|\partial_{z'}\tilde\phi_{u,v;\eps}(t,\cdot,1)|^2,
\frac{1+\eps^2v_{x}^2}{(u-v)^2}|\partial_{z'}\tilde\phi_{u,v;\eps}(t,\cdot,0)|^2\right)\label{defg}
\end{align}
and denote the components of $g_\eps$ by $g_{\eps,1}$ and $g_{\eps,2}$. Then clearly $g_\eps(\tilde u,\tilde v)(x)=g_\eps(u,v)(-x)$, $x\in I$, and \eqref{originalproblem3}--\eqref{originalproblem4} can be rewritten as
\begin{align}
u_t - u_{xx}&=-\lambda g_{\eps,1}(u,v), &x\in I,\, t>0, \label{transfeq3}\\
v_t - v_{xx}&=\mu g_{\eps,2}(u,v), &x\in I,\, t>0. \label{transfeq4}
\end{align}
Again we simplify notation by now omitting the index $\eps$. Given $(u,v)\in S_q(\kappa)$, we introduce a bounded linear operator $\mathcal A(u,v)\in\mathcal L(W^2_{2,D}(\Omega),L_2(\Omega))$ by setting
$$\mathcal A(u,v)w=-\widetilde{\Delta}_{u,v}w,\quad\forall w\in W^2_{2,D}(\Omega).$$
By \eqref{bddinv}, $\mathcal A(u,v)$ is invertible and its inverse $\mathcal A(u,v)^{-1}\in\mathcal L(L_2(\Omega),W^2_{2,D}(\Omega))$ satisfies
\beq\norm{\mathcal A(u,v)^{-1}}_{\mathcal L(L_2(\Omega),W^2_{2,D}(\Omega))}\leq c_3(\kappa,\eps).\label{boundA-1}\eeq
For nonzero $w\in W_{2,D}^2(\Omega)$ and $(u_1,v_1),(u_2,v_2)\in S_q(\kappa)$, we have, by \eqref{deltatilde} and the continuity of the mapping $L_q(\Omega)\cdot W_2^1(\Omega)\hookrightarrow L_2(\Omega)$,
\begin{align}
& \frac{\norm{\mathcal A(u_1,v_1)w-\mathcal A(u_2,v_2)w}_{L_2(\Omega)}}{\norm{w}_{W^2_{2,D}(\Omega)}}& \nonumber\\
& \leq 2\eps^2
\norm{\frac{z'(u_{1,x}-v_{1,x})+v_{1,x}}{u_1-v_1}-\frac{z'(u_{2,x}-v_{2,x})+v_{2,x}}{u_2-v_2}}_{L_\infty(\Omega)}\nonumber\\
& \quad + \norm{\frac{1+\eps^2[z'(u_{1,x}-v_{1,x})+v_{1,x}]^2}{(u_1-v_1)^2}
-\frac{1+\eps^2[z'(u_{2,x}-v_{2,x})+v_{2,x}]^2}{(u_2-v_2)^2}}_{L_\infty(\Omega)}\nonumber\\
& \quad + 2\eps^2\bigg\Vert\frac{(u_{1,x}-v_{1,x})[z'(u_{1,x}-v_{1,x})+v_{1,x}]}{(u_1-v_1)^2} \nonumber\\
& \hspace{2cm}-\frac{(u_{2,x}-v_{2,x})[z'(u_{2,x}-v_{2,x})+v_{2,x}]}{(u_2-v_2)^2}\bigg\Vert_{L_\infty(\Omega)} \nonumber\\
& \quad+\eps^2\norm{\frac{z'(u_{1,xx}-v_{1,xx})+v_{1,xx}}{u_1-v_1}
-\frac{z'(u_{2,xx}-v_{2,xx})+v_{2,xx}}{u_2-v_2}}_{L_q(\Omega)}\nonumber\\
& = 2\eps^2\norm{\alpha_1}_{L_\infty(\Omega)}+\norm{\alpha_2}_{L_\infty(\Omega)}+2\eps^2\norm{\alpha_3}_{L_\infty(\Omega)}+\eps^2\norm{\alpha_4}_{L_q(\Omega)}.\nonumber
\end{align}
Some elementary computations show that
\begin{align}
\alpha_1&=z'\left[\frac{u_{1,x}-u_{2,x}}{u_1-v_1}+u_{2,x}\frac{(u_2-u_1)-(v_2-v_1)}{(u_1-v_1)(u_2-v_2)}\right]\nonumber\\
&\hspace{1cm} +(1-z')\left[\frac{v_{1,x}-v_{2,x}}{u_1-v_1}
+v_{2,x}\frac{(u_2-u_1)-(v_2-v_1)}{(u_1-v_1)(u_2-v_2)}\right].\nonumber
\end{align}
We now make use of the continuity of the mappings $W^1_q(I)\cdot W^1_q(I)\hookrightarrow W^1_q(I)\hookrightarrow L_\infty(I)$ to see that $\norm{\alpha_1}_{L_\infty(\Omega)}$ can be bounded by a positive constant, depending only on $\kappa$, times $\norm{(u_1,v_1)-(u_2,v_2)}_{W_q^2(I)\times W_q^2(I)}$. The terms involving $\alpha_2,\ldots,\alpha_4$ can be treated similarly and we omit the tedious computations for the convenience of the reader. Finally, we get
\beq\label{absch1}\norm{\mathcal A(u_1,v_1)-\mathcal A(u_2,v_2)}_{\mathcal L(W^2_{2,D}(\Omega),L_2(\Omega))}\leq c_5(\kappa,\eps)\norm{(u_1,v_1)-(u_2,v_2)}_{W_q^2(I)\times W_q^2(I)}\eeq
and using the second resolvent identity, \eqref{boundA-1} and \eqref{absch1}, we also have
\begin{align}
&\norm{\mathcal A(u_1,v_1)^{-1}-\mathcal A(u_2,v_2)^{-1}}_{\mathcal L(L_2(\Omega),W^2_{2,D}(\Omega))}\nonumber\\
&\hspace{2cm}\leq c_6(\kappa,\eps)\norm{(u_1,v_1)-(u_2,v_2)}_{W_q^2(I)\times W_q^2(I)}.\label{absch2}
\end{align}
As before, one deduces
\begin{align}
\norm{f_{u_1,v_1}-f_{u_2,v_2}}_{L_2(\Omega)}&\leq 2\eps^2\norm{\alpha_3}_{L_2(\Omega)}+\eps^2\norm{\alpha_4}_{L_2(\Omega)}\nonumber\\
&\leq c_7(\kappa,\eps)\norm{(u_1,v_1)-(u_2,v_2)}_{W_q^2(I)\times W_q^2(I)}.\label{absch3}
\end{align}
From \eqref{absch3} and the fact that $(0,-1)\in\overline S_q(\kappa)$, we get $\norm{f_{u,v}}_{L_2(\Omega)}\leq 2c_7(\kappa,\eps)/\kappa$, for all $(u,v)\in S_q(\kappa)$. Now using \eqref{boundA-1}, \eqref{absch2} and \eqref{absch3}, we observe that
\begin{align}
\norm{\tilde\phi_{u_1,v_1}-\tilde\phi_{u_2,v_2}}_{W_2^2(\Omega)} & =\norm{\psi_{u_1,v_1}-\psi_{u_2,v_2}}_{W_2^2(\Omega)}\nonumber\\
& = \norm{\mathcal A(u_1,v_1)^{-1}f_{u_1,v_1}-\mathcal A(u_2,v_2)^{-1}f_{u_2,v_2}}_{W_{2,D}^2(\Omega)}\nonumber\\
& \leq c_8(\kappa,\eps)\norm{(u_1,v_1)-(u_2,v_2)}_{W_q^2(I)\times W_q^2(I)}\nonumber
\end{align}
and hence $S_q(\kappa)\to W_2^2(\Omega)\colon (u,v)\mapsto\tilde\phi_{u,v}$ is globally Lipschitz continuous. With the aid of \cite[Thm.~II-5.5]{N67} and the continuity of the pointwise multiplication $W_2^{1/2}(I)\cdot W_2^{1/2}(I)\hookrightarrow W_2^{2\sigma_1}(I)$, $2\sigma_1<1/2$, cf.\ \cite[Thm.~4.1]{Am91}, it follows that the mappings
$$S_q(\kappa)\to W_2^{2\sigma_1}(I),\quad (u,v)\mapsto\left|\partial_{z'}\tilde\phi_{u,v}(t,\cdot,1)\right|^2,\,(u,v)\mapsto\left|\partial_{z'}\tilde\phi_{u,v}(t,\cdot,0)\right|^2$$
are globally Lipschitz continuous. As $W_q^2(I)\hookrightarrow W_\infty^1(I)$, the mappings
$$S_q(\kappa)\to W_q^1(I),\quad (u,v)\mapsto\frac{1+\eps^2u_x^2}{(u-v)^2},\,(u,v)\mapsto\frac{1+\eps^2v_x^2}{(u-v)^2}$$
are globally Lipschitz continuous with a Lipschitz constant depending only on $\kappa$ and $\eps$. Finally the continuity of the pointwise multiplication $W_q^1(I)\cdot W_2^{2\sigma_1}(I)\hookrightarrow W_2^{2\sigma}(I)=W_{2,D}^{2\sigma}(I)$, $2\sigma<2\sigma_1<1/2$, cf.\ \cite[Thm.~4.1]{Am91}, implies that $g_\eps$ is globally Lipschitz continuous.

Note that the map $(u,v)\mapsto\tilde\phi_{u,v}\colon S_q(\kappa)\to W_2^2(\Omega)$ is analytic, since $\mathcal A\colon S_q(\kappa)\to\mathcal L(W^2_{2,D}(\Omega),L_2(\Omega))$ and hence $\mathcal A^{-1}\colon S_q(\kappa)\to\mathcal L(L_2(\Omega),W^2_{2,D}(\Omega))$ is analytic and by the analyticity of the right-hand side $(u,v)\mapsto f_{u,v}$, $S_q(\kappa)\to L_2(\Omega)$, of \eqref{phieqn1'}--\eqref{phieqn2'}. This immediately achieves that also $g_\eps$ is analytic. We have just proven the following proposition which is the analog of \cite[Prop.~2.1]{ELW12}.
\prop\label{propphig} Let $q\in(2,\infty)$, $\kappa\in(0,1/2)$ and $\eps>0$. For each $(u,v)\in S_q(\kappa)$ there is a unique solution $\tilde\phi_{u,v;\eps}\in W_2^2(\Omega)$ to the problem \eqref{phieqn1}--\eqref{phieqn2}. Moreover, with the definition $(\tilde u,\tilde v)(x)=(u,v)(-x)$, $x\in I$, we have that $\tilde\phi_{\tilde u,\tilde v;\eps}(t,x',z')=\tilde\phi_{u,v;\eps}(t,-x',z')$, $(x',z')\in\Omega$, $t>0$, and for $2\sigma\in[0,1/2)$, the mapping $g_\eps\colon S_q(\kappa)\to W_{2,D}^{2\sigma}(I)\times W_{2,D}^{2\sigma}(I)$ defined in \eqref{defg} is analytic, globally Lipschitz continuous and bounded with $g_\eps(0,-1)=(1,1)$.
\endprop\rm
Recall from \cite{mk13} that the Lipschitz continuity of the right-hand side of the equations on the free boundary was not needed for the stationary free boundary problem.

Now the boundary conditions \eqref{originalproblem5}--\eqref{originalproblem8} enter the game. For $p\in(1,\infty)$, we define a bounded linear operator $A_p\in\mathcal L(W_{p,D}^2(I),L_p(I))$ by setting $A_pw:=-w_{xx}$, for all $w\in W^2_{p,D}(I)$. As $A_r\subset A_p$, $r\geq p$, we simply write $A$ instead of $A_p$ in the following. Note that $-A$ is the generator of the heat semigroup $\{e^{-tA};\,t\geq 0\}$ on $L_p(I)$. In particular, $A$ is invertible.

With the definitions $\hat{v}=v+1$ and $\hat g_\eps(u,\hat v)=g_\eps(u,\hat v-1)=g_\eps(u,v)$, the equations \eqref{transfeq3} and \eqref{transfeq4} with the boundary conditions \eqref{originalproblem5}--\eqref{originalproblem8} read
\begin{align}
\left(\frac{d}{dt}+A\right)\begin{pmatrix}u\\\hat v\end{pmatrix}&=\begin{pmatrix}-\lambda & 0 \\ 0 & \mu\end{pmatrix}\hat g_\eps(u,\hat v),& x\in I,&\quad t>0,\label{evoleq1}\\
\begin{pmatrix}u\\\hat v\end{pmatrix}& = \begin{pmatrix}0\\0\end{pmatrix},& x\in\{1,-1\},&\quad t>0,\label{evoleq2}\\
\begin{pmatrix}u\\\hat v\end{pmatrix}& = \begin{pmatrix}u_0\\\hat{v}_0\end{pmatrix},& x\in I,&\quad t=0.\label{evoleq3}
\end{align}
Pick $u_0,v_0\in W^2_q(I)$, $q\in(2,\infty)$, such that $(u_0,v_0)(\pm 1)=(0,-1)$ and $-1\leq v_0<u_0\leq 0$ on $I$. Then there is $\kappa\in(0,1/4)$ such that $(u_0,v_0)\in\overline{S}_q(2\kappa)$. As explained in \cite[Lem.~2.3]{ELW12} there are $M\geq 1$ and $\omega>0$ so that
\beq\label{eq2.34}\norm{e^{-tA}}_{\mathcal L(W^2_{q,D}(I))}+t^{-\sigma+1+\tfrac{1}{2}(\tfrac{1}{2}-\tfrac{1}{q})}\norm{e^{-tA}}_{\mathcal L(W^{2\sigma}_{2,D}(I),W^2_{q,D}(I))}\leq Me^{-\omega t},\eeq
for $\tfrac{1}{2}-\tfrac{1}{q}<2\sigma<\tfrac{1}{2}$, $2\sigma\neq 1/q$. Let $\kappa_0=\kappa/M<\kappa$. Using the Lipschitz continuity of $\hat g_{\eps,i}$, for $i=1,2$, i.e.,
\begin{align}&\norm{\hat g_{\eps,i}(w_1,\hat w_2)-\hat g_{\eps,i}(w_3,\hat w_4)}_{W_{2,D}^{2\sigma}(I)}\nonumber\\
&\hspace{2cm}\leq c_9(\kappa,\eps)\norm{(w_1,\hat w_2)-(w_3,\hat w_4)}_{W_{q,D}^2(I)\times W_{q,D}^2(I)},
\label{gLip}
\end{align}
where $(w_1,w_2),(w_3,w_4)\in\overline S_q(\kappa_0)$, and that $(0,-1)\in\overline{S}_q(\kappa_0)$ and $g_\eps(0,-1)=(1,1)$, we obtain the bound
\beq\label{gbounded}\norm{\hat g_{\eps,i}(w_1,\hat w_2)}_{W_{2,D}^{2\sigma}(I)}\leq c_{10}(\kappa,\eps),\quad \forall(w_1,w_2)\in\overline{S}_q(\kappa_0),\eeq
for $i=1,2$. For $\tau>0$ we define the spaces $X_\tau:=C\left([0,\tau];\overline{S}_q(\kappa_0)+\{(0,1)\}\right)$ and for $t\in[0,\tau]$ and $(u,\hat v)\in X_\tau$ the map $F(u,\hat v)(t)=(F_1(u,\hat v)(t),F_2(u,\hat v)(t))^T$ given by
\begin{align}
F_1(u,\hat v)(t) & = e^{-tA}u_0-\lambda\int_0^te^{-(t-s)A}\hat g_{\eps,1}(u(s),\hat v(s))\,ds, \nonumber\\
F_2(u,\hat v)(t) & = e^{-tA}\hat v_0+\mu\int_0^te^{-(t-s)A}\hat g_{\eps,2}(u(s),\hat v(s))\,ds. \nonumber
\end{align}
We aim to apply the Contraction Mapping Theorem to the map $F$. Let
$$\mathcal I(\tau):=\int_0^\tau e^{-\omega s}s^{\sigma-1-\tfrac{1}{2}(\tfrac{1}{2}-\tfrac{1}{q})}\, ds.$$
Then $\mathcal I\to 0$ as $\tau\to 0$, $\mathcal I\to\mathcal I(\infty)<\infty$ for $\tau\to\infty$ and $\tau\mapsto\mathcal I(\tau)$ is increasing on $[0,\infty)$. Writing $m:=\max\{\lambda,\mu\}$ and using \eqref{eq2.34}--\eqref{gbounded}, we find
$$\norm{F_i(u,\hat v)(t)}_{W^2_{q,D}(I)}\leq\frac{M}{2\kappa}+mMc_{10}(\kappa,\eps)\mathcal I(\tau)$$
and
\begin{align}&\norm{F_i(u_1,\hat{v}_1)(t)-F_i(u_2,\hat{v}_2)(t)}_{W^2_{q,D}(I)}\nonumber\\
&\hspace{1.5cm}\leq mMc_{9}(\kappa,\eps)\mathcal I(\tau)\norm{(u_1,\hat{v}_1)-(u_2,\hat{v}_2)}_{C([0,\tau],W^2_{q,D}(I)\times W^2_{q,D}(I))},\nonumber\end{align}
for $i=1,2$, $(u_1,\hat v_1),(u_2,\hat v_2)\in X_\tau$ and $t\in[0,\tau]$. As $\norm{w}_{L_\infty(I)}\leq 2\norm{w}_{W^2_{q,D}(I)}$, for all $w\in W^2_{q,D}(I)$, $\hat g_{\eps,1},\hat g_{\eps,2}\geq 0$ and the heat semigroup is positivity preserving, we infer
\begin{align}
F_1(u,\hat v)(t) & \leq 0, \nonumber\\
F_2(u,\hat v)(t) & \geq 0, \nonumber\\
\hat F_1(u,\hat v)(t) - F_2(u,\hat v)(t) & \geq 4\kappa-4mMc_{10}(\kappa,\eps)\mathcal I(\tau). \nonumber
\end{align}
From this one concludes that there exists $\tau_0=\tau_0(\lambda,\mu,\kappa,\eps,q,\sigma)>0$ such that $F\colon X_{\tau_0}\to X_{\tau_0}$ is a contraction. It follows that there exists $T_\eps\in(\tau_0,\infty]$ and a unique maximal solution
$$\begin{pmatrix} u\\\hat v\end{pmatrix}=e^{-t A}\begin{pmatrix} u_0\\\hat v_0\end{pmatrix}+
\int_0^t e^{-(t-s)A}
\begin{pmatrix}-\lambda \hat g_{\eps,1}\\\mu \hat g_{\eps,2}\end{pmatrix}(u(s),\hat v(s))\,ds$$
to \eqref{evoleq1}--\eqref{evoleq3} on $[0,T_\eps)$ satisfying
$$u,\hat v\in C([0,T_\eps),W^2_{q,D}(I))\cap C((0,T_\eps),W_{2,D}^{2+2\sigma}(I))\cap C^1([0,T_{\eps}),L_q(I))$$
and
$$u(t,x)\leq 0,\quad\hat v(t,x) \geq 0,\quad \hat u(t,x)-\hat v(t,x)\geq 2\kappa_0,\quad (t,x)\in [0,T_\eps)\times I.$$
If, for any $\tau>0$, there is $\kappa(\tau)\in(0,1/2)$ and a solution $(u,\hat v)\in\overline S_q(\kappa(\tau))+\{(0,1)\}$ for $t\in[0,T_\eps)\cap[0,\tau]$, then $T_\eps=\infty$. Choosing $m$ suitably small, $m<m_1$, where $m_1=m_1(\kappa,\eps,q,\sigma)>0$, we obtain
$$mM\max\{c_{9},c_{10}\}I(\infty)<1<\frac{1}{2\kappa_0}\quad\text{and}\quad 2mMc_{10}I(\infty)\leq\kappa_0,$$
so that the map $F\colon X_\tau\to X_\tau$ is a contraction for any $\tau>0$. In particular, there exists a unique global solution $(u,\hat v)\in\overline S_q(\kappa_0)+\{(0,1)\}$. Finally, Proposition~\ref{propphig} and uniqueness of the solution imply that, for $u_0,\hat{v}_0$ even, the associated solution $(u,\hat v)$ to \eqref{evoleq1}--\eqref{evoleq3} is even on $[0,T_\eps)\times I$. Up to the transformation $v=\hat v-1$ and up to pulling the solution $(u,v,\tilde\phi)$ back to $\Omega_{u(t),v(t)}$, this completes the proof of the following theorem.
\thm\label{thm1} Let $q\in(2,\infty)$, $\eps>0$ and initial values $u_0,v_0\in W^2_{q}(I)$ with $-1\leq v_0<u_0\leq 0$ on $I$ and $(u_0,v_0)(\pm 1)=(0,-1)$ be given. Then:
\begin{enumerate}
\item[(i)] For any $\lambda,\mu>0$, there is a unique maximal solution $(u_\eps,v_\eps,\phi_\eps)$ to \eqref{originalproblem1}--\eqref{originalproblem8} with regularity
$$u_\eps,v_\eps\in C([0,T_\eps),W^2_{q}(I))\cap C^1([0,T_\eps),L_q(I)),\,\phi\in W_2^2(\Omega_{u_\eps(t),v_\eps(t)})$$
so that $-1\leq v_\eps<u_\eps\leq 0$ on $[0,T_\eps)\times I$ and $T_\eps>0$ is maximal.
\item[(ii)] If for each $\tau>0$ there is $\kappa(\tau)\in(0,1/2)$ such that $(u_\eps(t),v_\eps(t))\in S_q(\kappa(\tau))$ for $t\in[0,T_\eps)\cap[0,\tau]$, then the solution exists globally in time, i.e., $T_\eps=\infty$.
\item[(iii)] If $u_0$ and $v_0$ are even functions on $I$, then $(u_\eps,v_\eps,\phi_\eps)$ is even in $x$ on $[0,T_\eps)\times I$.
\item[(iv)] Given $\kappa\in(0,1/2)$ and $(u_0,v_0)\in S_q(\kappa)$, there exist $m_1=m_1(\kappa,\eps)>0$ and $\kappa_0=\kappa_0(\kappa,\eps)>0$ such that, for $\lambda,\mu<m_1$, $T_\eps=\infty$ and $(u_\eps(t),v_\eps(t))\in S_q(\kappa_0)$ for all $t\geq 0$.
\end{enumerate}
\endthm\rm
Before we proceed, we prepare the following lemma about some elementary properties of the solution $(u,v,\phi)$ to \eqref{originalproblem1}--\eqref{originalproblem8}. Theorem~\ref{thm1} and Lemma~\ref{lemxxx} are the analogs of Theorem~1.1, Theorem~1.2(i) and Proposition~2.4 of \cite{ELW12}.
\lem\label{lemxxx} Let $q\in(2,\infty)$, $\eps,\lambda,\mu>0$ and initial values $u_0,v_0\in W_q^2(I)$ with $(u_0,v_0)(\pm 1)=(0,-1)$ and $-1\leq v_0<u_0\leq 0$ on $I$ be given. Let $(u_\eps,v_\eps,\phi_\eps)$ denote the associated maximal solution of \eqref{originalproblem1}--\eqref{originalproblem8} satisfying the properties stated in Theorem~\ref{thm1}. Then, for all $(t,x,z)\in[0,T_\eps)\times\overline\Omega_{u,v}$,
\beq\label{boundphi}-1\leq\phi_\eps(t,x,z)\leq 1\eeq
and
\begin{align}
(\partial_x\phi_\eps)(t,x,u(t,x))&=-(\partial_xu_\eps)(t,x)(\partial_z\phi_\eps)(t,x,u(t,x)),\label{chainrulex}\\
(\partial_x\phi_\eps)(t,x,v(t,x))&=-(\partial_xv_\eps)(t,x)(\partial_z\phi_\eps)(t,x,v(t,x)).\label{chainrulez}
\end{align}
\endlem\rm
\proof The bounds \eqref{boundphi} are obtained from the maximum principle applied to the constant functions $\pm 1$ and the function $\phi_\eps$. Finally, differentiating the equations $\phi_\eps(t,x,u_\eps(t,x))=1$ and $\phi_\eps(t,x,v_\eps(t,x))=0$ with respect to $x$ and applying the chain rule, we immediately obtain \eqref{chainrulex} and \eqref{chainrulez}.
\endproof
\section{Non-existence of global solutions}\label{sec_nonex}
%
In Section~\ref{sec_lwp} we have proved the local existence of solutions to \eqref{originalproblem1}--\eqref{originalproblem8}. Let us now discuss criteria for the non-existence of global solutions. Let $(u,v,\phi)$ denote the maximal solution of \eqref{originalproblem1}--\eqref{originalproblem8} with initial values $u_0,v_0$ satisfying the properties stated in Theorem~\ref{thm1}; to simplify notation, we omit the index $\eps$ in this section again. Here, we show that there is a critical value $m_2(\eps)>0$ such that for $\max\{\lambda,\mu\}>m_2$, the maximal existence time $T_\eps>0$ of the solution $(u,v,\phi)$ is finite.
\thm\label{thmnonex} Let $q\in(2,\infty)$, $\eps,\lambda,\mu>0$ and initial values $u_0,v_0\in W_q^2(I)$ with $(u_0,v_0)(\pm 1)=(0,-1)$ and $-1\leq v_0<u_0\leq 0$ on $I$ be given. Let $(u,v,\phi)$ denote the associated maximal solution of \eqref{originalproblem1}--\eqref{originalproblem8} with initial values $u_0,v_0$ and maximal existence time $T_\eps>0$ according to Theorem~\ref{thm1}. Then for
$$\max\{\lambda,\mu\}>m_2:=\pi^4(1+\eps^2)^2$$
we have that $T_\eps<\infty$.
\endthm\rm
\proof Let $\zeta_1\colon\overline I\to[0,\pi/4]$, $\zeta_1(x):=\tfrac{\pi}{4}\cos(\tfrac{\pi}{2}x)$ and $\mu_1:=\tfrac{\pi^2}{4}$ so that
$$-\partial_x^2\zeta_1= \mu_1\zeta_1,\;x\in I,\quad \zeta_1(\pm 1)=0,\quad\norm{\zeta_1}_{L_1(I)}=1,$$
i.e., $\mu_1$ is the principal eigenvalue of $-\partial_x^2$ acting on $L_2(I)$. For some $\alpha\in(0,1)$ and $t\in[0,T_\eps)$, let
$$E_{\alpha}(t):=\int_I\zeta_1(x)(u+\tfrac{\alpha}{2}u^2)(t,x)\, dx.$$
In \cite[Sec.~3]{ELW12}, it is shown for the problem with one free boundary that
\beq\label{derivneg}\frac{d}{dt}E_{\alpha}\leq\mu_1 + \frac{4\alpha\beta}{\eps^2p}\left[\frac{1}{p}\mu_1\eps^2+\frac{1}{4\beta}p-\frac{1}{1+E_\alpha}\right],\eeq
where $\alpha=\eps^2/(1+\eps^2)$, $\beta=\sqrt\lambda/2$ and $p=1+2\mu_1\eps^2$, and it is proven that, for $\lambda>m_2$, the right-hand side of the inequality \eqref{derivneg} can be estimated by a negative constant, so that finiteness of $T_\eps$ follows immediately by integrating \eqref{derivneg} over $[0,T_\eps)$. Using Lemma~\ref{lemxxx}, it is straightforward to generalize these arguments for $E_\alpha(t)$ in the problem with two free boundaries, with the obvious changes, so that \eqref{derivneg} also holds true for the problem \eqref{originalproblem1}--\eqref{originalproblem8}. In the case $\mu>\lambda$, we make use of the following symmetry of the problem \eqref{originalproblem1}--\eqref{originalproblem8}: If $(u,v,\phi)$ is a solution to \eqref{originalproblem1}--\eqref{originalproblem8} with the parameters $(\lambda,\mu)$, then $(U,V,\varphi)$ defined by
$$U=-v-1,\quad V=-u-1,\quad \varphi(t,x,z)=1-\phi(t,x,-z-1)$$
is a solution to \eqref{originalproblem1}--\eqref{originalproblem8} with the parameters $(\mu,\lambda)$. In view of this symmetry and the inequality \eqref{derivneg}, the proof is completed.
\endproof
Theorem~\ref{thmnonex} shows that, for $\lambda$ or $\mu$ sufficiently large, the problem \eqref{originalproblem1}--\eqref{originalproblem8} cannot have a stationary solution. Applying a technique presented in \cite[Thm.~3]{WL12}, we obtain a more precise value of the threshold for the parameters $\lambda$ and $\mu$ to guarantee the non-existence of stationary solutions of \eqref{originalproblem1}--\eqref{originalproblem8}. Recall from \cite{mk13} that the stationary problem \eqref{originalproblem1}--\eqref{originalproblem8} possesses even solutions in $W^2_\infty(I)^2\times W^2_2(\Omega_{u,v})$.
\thm Let $\eps>0$. There exists $\xi_0(\eps)\in(0,\tfrac{\pi}{2\eps})$ such that for $\max\{\lambda,\mu\}>\xi_0(\eps)$ the stationary problem \eqref{originalproblem1}--\eqref{originalproblem6} possesses no even solution $(u,v,\phi)$ of regularity $u,v\in W^2_\infty(I)$ and $\phi\in W_2^2(\Omega_{u,v})$ such that $-1\leq v<u\leq 0$ on $I$. In addition $\xi_0(\eps)\to 2$ for $\eps\to 0$.
\endthm\rm
\proof
It follows from \eqref{originalproblem3} and \eqref{originalproblem4} that, for all $x\in I$,
\begin{align}
u_{xx}(x)&=\lambda(1+\eps^2|u_{x}(x)|^2)|\phi_z(x,u(x))|^2, \label{stationaryuxx}\\
v_{xx}(x)&=-\mu(1+\eps^2|v_{x}(x)|^2)|\phi_z(x,v(x))|^2. \label{stationaryvxx}
\end{align}
For reasons of convexity, cf.\ \cite{mk13}, $S_1(x,z):=1+z-u(x)$ is a supersolution and $S_2(x,z):=z-v(x)$ is a subsolution for the elliptic operator $-\Delta_\eps$ satisfying
\begin{align}
S_1(\pm 1,z)&=1+z&\hspace{-1.6cm}=&\;\phi(\pm 1,z),\nonumber\\
S_1(x,u(x))&=1&\hspace{-1.6cm}=&\;\phi(x,u(x)),\nonumber\\
S_1(x,v(x))&=1+v(x)-u(x)\geq 0&\hspace{-1.6cm}=&\;\phi(x,v(x)),\nonumber\\
-\Delta_\eps S_1(x,z)&=\eps^2 u_{xx}\geq 0&\hspace{-1.6cm}=&\;-\Delta_\eps\phi\nonumber
\end{align}
and
\begin{align}
S_2(\pm 1,z)&=1+z&\hspace{-1.6cm}=&\;\phi(\pm 1,z),\nonumber\\
S_2(x,u(x))&=u(x)-v(x)\leq 1&\hspace{-1.6cm}= &\;\phi(x,u(x)),\nonumber\\
S_2(x,v(x))&=0&\hspace{-1.6cm}=&\;\phi(x,v(x)),\nonumber\\
-\Delta_\eps S_2(x,z)&=\eps^2 v_{xx}\leq 0&\hspace{-1.6cm}=&\;-\Delta_\eps\phi.\nonumber
\end{align}
The weak maximum principle implies that, for all $(x,z)\in\overline\Omega_{u,v}$,
$$\phi(x,z)-S_1(x,z)\leq\max_{\overline\Omega_{u,v}}\left\{\phi-S_1\right\}=\max_{\partial\Omega_{u,v}}\left\{\phi-S_1\right\}\leq 0$$
and that
$$S_2(x,z)-\phi(x,z)\leq\max_{\overline\Omega_{u,v}}\left\{S_2-\phi\right\}=\max_{\partial\Omega_{u,v}}\left\{S_2-\phi\right\}\leq 0.$$
Hence, for fixed $x\in I$ and $z\in(v(x),u(x))$,
$$\frac{\phi(x,z)-\phi(x,u(x))}{z-u(x)}\geq 1\quad\text{and}\quad\frac{\phi(x,z)-\phi(x,v(x))}{z-v(x)}\geq 1,$$
and sending $z$ to $u(x)$ and $v(x)$ respectively, we conclude that $\phi_z(x,u(x))\geq 1$ and $\phi_z(x,v(x))\geq 1$, for all $x\in I$.
Then \eqref{stationaryuxx} and \eqref{stationaryvxx} imply
\begin{align}
u_{xx}(x)&\geq\lambda(1+\eps^2|u_{x}(x)|^2), \label{stationaryuxx2}\\
v_{xx}(x)&\leq-\mu(1+\eps^2|v_{x}(x)|^2). \label{stationaryvxx2}
\end{align}
Let
$$\Lambda_\eps(\xi):=1+\frac{1}{\eps^2\xi}\ln(\cos(\eps\xi)),\quad\xi\in\left(0,\frac{\pi}{2\eps}\right).$$
As explained in the proof of \cite[Thm.~3]{WL12}, $\Lambda_\eps$ possesses a unique zero $\xi_0(\eps)$ such that $\Lambda_\eps(\xi)<0$ for $\xi_0(\eps)<\xi<\frac{\pi}{2\eps}$ and $\xi_0(\eps)\to 2$ for $\eps\to 0$. Moreover, it has been shown that integrating \eqref{stationaryuxx2} twice leads to
$$u(0)\leq\frac{1}{\eps^2\lambda}\ln(\cos(\eps\lambda x)),\quad x\in\left[0,\min\left\{1,\frac{\pi}{2\eps\lambda}\right\}\right)$$
so that, for $\lambda\geq\tfrac{\pi}{2\eps}$ and $x\to\tfrac{\pi}{2\eps\lambda}$, $u(0)=-\infty$, and, for $x\to 1$, $u(0)\leq\Lambda_\eps(\lambda)-1\leq-1$ on $[\xi_0(\eps),\tfrac{\pi}{2\eps})$, which are both contradictions.

%
Clearly, any $C^1$-smooth even function on $I$ has vanishing derivative at $x=0$. Then $v_x(0)=0$ and integrating \eqref{stationaryvxx2} over $[0,x]$ yield
$$\arctan(\eps v_x(x))\leq -\eps\mu x,\quad x\in[0,1),$$
or equivalently
$$v_x(x)\leq-\frac{1}{\eps}\tan\left(\eps\mu x\right),\quad x\in\left[0,\min\left\{1,\frac{\pi}{2\eps\mu}\right\}\right).$$
Integrating once more and using that $v(0)<0$, we arrive at
$$v(x)<\frac{1}{\eps^2\mu}\ln(\cos(\eps\mu x)),\quad x\in\left[0,\min\left\{1,\frac{\pi}{2\eps\mu}\right\}\right).$$
Assuming $\mu\geq\frac{\pi}{2\eps}$ and sending $x\to\frac{\pi}{2\eps\mu}$, we see that $v(\frac{\pi}{2\eps\mu})=-\infty$ which is clearly contradicting $v\geq-1$. Assuming $\mu<\frac{\pi}{2\eps}$ and sending $x\to 1$, we get $v(1)\leq\Lambda_{\eps}(\mu)-1$ and for $\mu>\xi_0(\eps)$, we conclude $-1=v(1)<-1$, which is again a contradiction. This completes the proof of our theorem.\endproof
\section{Asymptotic stability}\label{sec_stab}
Fix $q\in(2,\infty)$, $\kappa\in(0,1/2)$, $\eps>0$ and $2\sigma\in(\frac{1}{2}-\frac{1}{q},\frac{1}{2})$. Recall that $g_\eps$ defined in \eqref{defg} is an analytic map $S_q(\kappa)\to W_{2,D}^{2\sigma}(I)\times W_{2,D}^{2\sigma}(I)\hookrightarrow L_q(I)\times L_q(I)$. Moreover the operator $A=-\partial_x^2$ with $D(A)=W_{q,D}^2(I)$ is invertible and $-A$ generates the heat semigroup on $L_q(I)$.

Define $F\colon\R^2\times S_q(\kappa)\to W_{q,D}^2(I)\times W^2_{q}(I)$ by setting
$$
F(\Lambda,U)=\begin{pmatrix}U_1\\U_2\end{pmatrix} + \begin{pmatrix}\Lambda_1 & 0 \\ 0 & -\Lambda_2\end{pmatrix}A^{-1}g_\eps(U_1,U_2).
$$
Then $F(0,0)=(0,0)$ and, for all $W\in S_q(\kappa)$,
$$[D_UF(0,0)]W = \lim_{t\to 0}\frac{1}{t}(F(0,tW)-F(0,0)) = W.$$
According to the Implicit Function Theorem, there is $\delta>0$ and an analytic function $[0,\delta)^2\to S_q(\kappa)$, $\Lambda\mapsto U_\Lambda$, such that $F(\Lambda,U_\Lambda)=0$. For $\Lambda\neq (0,0)$, let $\Phi_\Lambda$ denote the associated potential solving \eqref{originalproblem1} and \eqref{originalproblem2}
with $u,v$ replaced by $U_{\Lambda,1},U_{\Lambda,2}$. Writing $\Lambda=(\lambda,\mu)$, $(U_\Lambda,\Phi_\Lambda)$ is a stationary solution of \eqref{originalproblem1}--\eqref{originalproblem6} as $F(\Lambda,U_\Lambda)=0$ and $U_{\Lambda,t}=0$ imply the equations corresponding to \eqref{originalproblem3} and \eqref{originalproblem4} and also \eqref{originalproblem5} and \eqref{originalproblem6} are satisfied. With the notation $U=(u,v)$ and $\hat U=(u,\hat v)$, equations \eqref{originalproblem3} and \eqref{originalproblem4} read
$$\hat U_t+A\hat U=\begin{pmatrix}-\Lambda_1 & 0 \\ 0 & \Lambda_2\end{pmatrix}g_\eps(U).$$
Setting $\hat V=U-U_\Lambda=\hat U-\hat U_\Lambda$, $\Lambda\in(0,\delta)^2$, and
$$B_\Lambda:=-\begin{pmatrix}-\Lambda_1 & 0 \\ 0 & \Lambda_2\end{pmatrix}Dg_\eps(U_\Lambda)\in\mathcal L(W^2_{q,D}(I)\times W^2_{q,D}(I),L_q(I)\times L_q(I)),$$
we obtain the linearization
\beq\label{linearization}\hat V_t+(A+B_\Lambda)\hat V=\begin{pmatrix}-\Lambda_1 & 0 \\ 0 & \Lambda_2\end{pmatrix}
\left(g_\eps(\hat V+U_\Lambda)-g_\eps(U_\Lambda)-Dg_\eps(U_\Lambda)\hat V\right),\eeq
and, denoting the right-hand side of \eqref{linearization} by $G_\Lambda(\hat V)$, the initial value problem
\begin{align}
\hat V_t+(A+B_\Lambda)\hat V&=G_\Lambda(\hat V),\quad t>0,\nonumber\\
\hat V(0)&=\hat V_0,\nonumber
\end{align}
where $G_\Lambda\in C^2(\mathcal O_\Lambda,L_q(I)\times L_q(I))$ is defined on an open zero neighborhood $\mathcal O_\Lambda\subset W_{q,D}^2(I)\times W_{q,D}^2(I)$ such that $U_\Lambda+\mathcal O_\Lambda\subset S_q(\kappa)$. Moreover $G_\Lambda(0)=0$ and $DG_\Lambda(0)=0$. It follows from a line of arguments similar to what is presented in Section 4 of \cite{ELW12} that
$$\lim_{\Lambda\to 0}\norm{B_\Lambda}_{\mathcal L(W^2_{q,D}(I)\times W^2_{q,D}(I),L_q(I)\times L_q(I))}=0$$
implies that $-(A+B_\Lambda)$ generates an analytic semigroup on $L_q(I)\times L_q(I)$ with a negative spectral bound. Now the following theorem is an immediate consequence of \cite[Thm.\ 9.1.2]{Lunardi}. 
\thm\label{thm_stability} Let $q\in(2,\infty)$, $\kappa\in(0,1/2)$ and $\eps>0$ be fixed.
\begin{enumerate}
\item[(i)] There are $\delta(\kappa)>0$ and an analytic function $[0,\delta)^2\to W^2_{q,D}(I)\times W_{q}^2(I)$, $\Lambda\to U_\Lambda=(U_{\Lambda,1},U_{\Lambda,2})$, such that, for each $\Lambda=(\lambda,\mu)\in(0,\delta)^2$, $(U_\Lambda,\Phi_\Lambda)$ is the unique steady state of \eqref{originalproblem1}--\eqref{originalproblem6} with $U_\Lambda\in S_q(\kappa)$ and $\Phi_\Lambda\in W_2^2(\Omega_{U_{\Lambda,1},U_{\Lambda,2}})$. Moreover, $U_{\Lambda,1}$ and $-U_{\Lambda,2}$ are convex and even for all $\Lambda\in(0,\delta)^2$ and $U_{(0,0)}=(0,0)$.
\item[(ii)] Let $\Lambda\in(0,\delta)^2$. There are $\omega_0,r,R>0$ such that for each pair of initial values $u_0,v_0\in W_{q}^2(I)$ satisfying $(u_0,v_0)(\pm 1)=(0,-1)$, $-1\leq v_0<u_0\leq 0$ and $\norm{(u_0,v_0)-U_{\Lambda}}_{W_{q,D}^2(I)\times W_{q,D}^2(I)}<r$,
    the associated solution $(u,v,\phi)$ to \eqref{originalproblem1}--\eqref{originalproblem8} exists globally in time and
    \begin{align}
    &\norm{(u,v)-U_\Lambda}_{W_{q,D}^2(I)\times W_{q,D}^2(I)}+\norm{(u_t,v_t)}_{L_q(I)\times L_q(I)}\nonumber\\
    &\hspace{2cm}\leq Re^{-\omega_0t}\norm{(u_0,v_0)-U_{\Lambda}}_{W_{q,D}^2(I)\times W_{q,D}^2(I)},\quad\forall t\geq 0.\nonumber
    \end{align}
\end{enumerate}
\endthm\rm
\rem As a consequence of the above theorem and the Lipschitz continuity of $(u,v)\mapsto\tilde\phi_{u,v}$, we also have, under the assumptions of Theorem \ref{thm_stability}, that $\tilde\phi_{u,v}$ converges exponentially to $\tilde\phi_{U_{\Lambda,1},U_{\Lambda,2}}$ as $t\to\infty$, i.e.,
$$\norm{\tilde\phi_{u,v}-\tilde\phi_{U_{\Lambda,1},U_{\Lambda,2}}}_{W_2^2(\Omega)}\leq R'e^{-\omega_0t}\norm{(u_0,v_0)-U_{\Lambda}}_{W_{q,D}^2(I)\times W_{q,D}^2(I)},\quad\forall t\geq 0,$$
with a positive constant $R'$.
\endrem\rm
\section{The small aspect ratio limit}\label{sec_sar}
In this section, we examine the connection between the original problem \eqref{originalproblem1}--\eqref{originalproblem8} and the vanishing aspect ratio model \eqref{sarm1}--\eqref{sarm7}. Let $\lambda,\mu>0$, $q\in(2,\infty)$ and $\kappa\in(0,1/2)$ be fixed. With $M$ in \eqref{eq2.34}, define $\kappa_1:=\kappa/(2M)<\kappa$. For $\eps>0$, let $(u_\eps,v_\eps,\phi_\eps)(t)$ denote the solution of \eqref{originalproblem1}--\eqref{originalproblem8} on $[0,T_\eps)$, for $(u_0,v_0)\in S_q(\kappa)$ with $u_0\leq 0$ and $v_0\geq-1$ given, cf.\ Theorem~\ref{thm1} and its proof. As the solution is continuous in time,
$$\tau_\eps:=\sup\left\{t\in[0,T_\eps);(u_\eps(s),v_\eps(s))\in\overline S_q(\kappa_1);\,\forall s\in[0,t]\right\}$$
is positive. Moreover, $T_\eps\geq\tau_\eps$. We then have
$$u_\eps(t)-v_\eps(t)\geq 2\kappa_1,\quad -1\leq v_\eps(t)<u_\eps(t)\leq 0\quad\text{on }[0,\tau_\eps]\times[-1,1],$$
and, by the continuous embedding $W_q^2(I)\hookrightarrow W_\infty^1(I)$,
$$\norm{u_\eps(t)}_{W_q^2(I)} + \norm{v_\eps(t)}_{W_q^2(I)} + \norm{u_\eps(t)}_{W_\infty^1(I)} + \norm{v_\eps(t)}_{W_\infty^1(I)}
\leq C_1,\quad\forall t\in[0,\tau_\eps].$$
Again, we denote by $C_1,C_2,C_3,\ldots$ a sequence of positive constants. Henceforth, we choose $\eps$ sufficiently small, i.e., $\eps$ smaller than some $\eps_1>0$, so that
\beq\label{eq5.4}\eps_1^2\left(\norm{u_{\eps,x}(t)}_{L_\infty(I)}+2\norm{v_{\eps,x}(t)}_{L_\infty(I)}\right)^2\leq\frac{1}{2},
\quad\forall(t,\eps)\in[0,\tau_\eps]\times(0,\eps_1].\eeq
For $(t,x',z')\in[0,\tau_\eps]\times\overline\Omega$, we recall the definition $\psi_\eps(t,x',z')=\tilde\phi_\eps(t,x',z')-z'$, where
$\tilde\phi_\eps(t) = \theta_*(u(t),v(t))\phi_\eps(t) = \phi_\eps(t)\circ T_{u_\eps(t),v_\eps(t)}^{-1}$ with the transformation $T_{u_\eps(t),v_\eps(t)}$ in \eqref{defT}. Also recall that, by Lemma~\ref{lemxxx}, $-1\leq\tilde\phi_\eps\leq 1$, so that $-2\leq\psi_\eps\leq 1$ on $[0,\tau_\eps]\times\overline\Omega$. The function $f_\eps(t,x',z')=\widetilde\Delta_{u_\eps(t),v_\eps(t);\eps} z'$ has been computed in \autoref{deff}. To simplify notation, we will write $(x,z)$ for points in $\Omega$ henceforth, since we do not need to distinguish between points in $\Omega$ and $\Omega_{u,v}$ here.

In what follows, we need control of the $L_2(\Omega)$-norm of $f_\eps$ and therefore we prepare the estimates
\begin{align}
\norm{f_\eps}_{L_q(\Omega)} & \leq \frac{2\eps^2}{4\kappa_1^2}\norm{z(u_{\eps,x}-v_{\eps,x})+v_{\eps,x}}_{L_\infty(\Omega)}\norm{u_{\eps,x}-v_{\eps,x}}_{L_q(I)} \nonumber\\
& \hspace{1cm} + \frac{\eps^2}{2\kappa_1}\left(\norm{u_{\eps,xx}}_{L_q(I)}+2\norm{v_{\eps,xx}}_{L_q(I)}\right)\leq C_2\eps^2\nonumber
\end{align}
and, with the aid of H\"older's inequality,
$$\norm{f_\eps}_{L_p(\Omega)}\leq 2^{\frac{q-p}{qp}}\norm{f_\eps}_{L_q(\Omega)}\leq C_3\eps^2,\quad\forall p\in[1,q].$$
Our next lemma provides some important bounds on the function $\psi_\eps$ and its derivatives. Compared to the analysis of the stationary case, it can be seen as a generalization of \cite[Lem.~7]{mk13} and \cite[Lem.~5.1]{ELW12}.
\lem\label{lemyyy} There exist positive constants $\eps_*$ and $K_1$ such that, for all $\eps<\eps_*$ and $t\in[0,\tau_\eps]$,
\begin{align}
\norm{\partial_{x}\psi_\eps(t)}_{L_2(\Omega)}+\frac{1}{\eps}\left(\norm{\psi_\eps(t)}_{L_2(\Omega)}+
\norm{\partial_{z}\psi_\eps(t)}_{L_2(\Omega)}\right)&\leq K_1,\label{eq5.5}\\
\frac{1}{\eps}\norm{\partial_{x}\partial_{z}\psi_\eps(t)}_{L_2(\Omega)}+\frac{1}{\eps^2}\norm{\partial_{z}^2\psi_\eps(t)}_{L_2(\Omega)}&\leq K_1,\label{eq5.6}\\
\frac{1}{\eps}\left(\norm{\partial_{z}\psi_\eps(t,\cdot,0)}_{W^{1/2}_2(I)}+\norm{\partial_{z}\psi_\eps(t,\cdot,1)}_{W^{1/2}_2(I)}\right)&\leq K_1.\label{eq5.7}
\end{align}
\endlem\rm
\proof As $-\widetilde\Delta_\eps\psi_\eps=f_\eps$ and $\psi_\eps|_{\partial\Omega}=0$, we can use the divergence form of $-\widetilde\Delta_\eps$ and integration by parts to obtain
\begin{align}
\int_\Omega f_\eps\psi_\eps\dx\dz & = \eps^2\int_\Omega\left(\partial_{x}\psi_\eps-\frac{z(u_{\eps,x}-v_{\eps,x})+v_{\eps,x}}{u_{\eps}-v_{\eps}}\partial_{z}\psi_\eps\right)^2dx\dz \nonumber\\
& \quad + \int_\Omega\frac{|\partial_{z}\psi_\eps|^2}{(u_{\eps}-v_{\eps})^2}\dx\dz\nonumber\\
& \quad + \eps^2 \int_\Omega(u_{\eps,x}-v_{\eps,x})\frac{z(u_{\eps,x}-v_{\eps,x})+v_{\eps,x}}{(u_\eps-v_\eps)^2}(\partial_{z}\psi_\eps)\psi_\eps\dx\dz\nonumber\\
& \quad -\eps^2 \int_\Omega\frac{u_{\eps,x}-v_{\eps,x}}{u_\eps-v_\eps}(\partial_x\psi_\eps)\psi_\eps\dx\dz.\nonumber
\end{align}
Now using the inequality $(r-s)^2\geq\frac{r^2}{2}-s^2$ we get
\begin{align}
& \int_\Omega f_\eps\psi_\eps\dx\dz \geq \frac{\eps^2}{2}\norm{\partial_{x}\psi_\eps}^2_{L_2(\Omega)}\nonumber\\
& \quad + \left(1- \eps^2\norm{z(u_{\eps,x}-v_{\eps,x})+v_{\eps,x}}^2_{L_\infty(\Omega)}\right)\int_\Omega\frac{|\partial_{z}\psi_\eps|^2}{(u_\eps-v_\eps)^2}\dx\dz \nonumber\\
& \quad - \frac{\eps^2}{2}\norm{z(u_{\eps,x}-v_{\eps,x})^2+v_{\eps,x}(u_{\eps,x}-v_{\eps,x})}_{L_\infty(\Omega)}
\int_\Omega\left(\frac{|\partial_{z}\psi_\eps|^2+|\psi_\eps|^2}{(u_\eps-v_\eps)^2}\right)dx\dz\nonumber\\
& \quad -\frac{\eps^2}{4}\int_\Omega|\partial_x\psi_\eps|^2\dx\dz-\eps^2
\norm{\frac{u_{\eps,x}-v_{\eps,x}}{u_\eps-v_\eps}}_{L_\infty(I)}^2\int_\Omega|\psi_\eps|^2\dx\dz\nonumber\\
& \geq \frac{\eps^2}{4}\norm{\partial_{x}\psi_\eps}_{L_2(\Omega)}^2+(1-C_4\eps^2)\norm{\partial_{z}\psi_\eps}_{L_2(\Omega)}^2-C_5\eps^2.\nonumber
\end{align}
Since $$\int_\Omega f_\eps\psi_\eps\dx\dz\leq\norm{f_\eps}_{L_2(\Omega)}\norm{\psi_\eps}_{L_2(\Omega)}\leq C_6\eps^2,$$ we have, for sufficiently small $\eps$, i.e., $\eps$ smaller than some $\eps_2>0$,
\beq\label{eq5.10}C_7\eps^2\geq \eps^2\norm{\partial_{x}\psi_\eps}_{L_2(\Omega)}^2+\norm{\partial_{z}\psi_\eps}_{L_2(\Omega)}^2,\eeq
and \eqref{eq5.10} shows that $\norm{\partial_{x}\psi_\eps}_{L_2(\Omega)}$ and $\frac{1}{\eps}\norm{\partial_{z}\psi_\eps}_{L_2(\Omega)}$ are bounded by a positive constant. Since $\norm{\psi_\eps}_{L_2(\Omega)}\leq\norm{\partial_{z}\psi_\eps}_{L_2(\Omega)}$, it is clear that $\frac{1}{\eps}\norm{\psi_\eps}_{L_2(\Omega)}$ is bounded by a positive constant. This implies \eqref{eq5.5}. Setting $\zeta_\eps:=\partial_{z}^2\psi_\eps$ and $\omega_\eps:=\partial_{x}\partial_{z}\psi_\eps$, it follows from integrating the equation $-\widetilde\Delta_\eps{\psi_\eps}=f_\eps$, as explained in \cite{mk13}, that
\begin{align}
&\int_\Omega f_\eps(1-\partial_{z}\psi_\eps)\zeta_\eps\dx\dz \nonumber\\
& \qquad = \eps^2\int_\Omega\left(\omega_\eps-\frac{z(u_{\eps,x}-v_{\eps,x})+v_{\eps,x}}{u_\eps-v_\eps}\zeta_\eps\right)^2dx\dz
+\int_\Omega\frac{\zeta_\eps^2}{(u_\eps-v_\eps)^2}\dx\dz.\nonumber
\end{align}
Using again the inequality $(r-s)^2\geq\frac{r^2}{2}-s^2$ and \eqref{eq5.4} it follows that
\begin{align}
&\int_\Omega f_\eps(1-\partial_{z}\psi_\eps)\zeta_\eps\dx\dz \nonumber\\
&\quad\geq \int_\Omega\left[\frac{\zeta_\eps^2}{(u_\eps-v_\eps)^2}+\frac{\eps^2}{2}\omega_\eps^2-\eps^2\zeta_\eps^2
\left(\frac{z(u_{\eps,x}-v_{\eps,x})+v_{\eps,x}}{u_\eps-v_\eps}\right)^2\right]dx\dz\nonumber\\
&\quad\geq \frac{1}{2}\left(\norm{\zeta_\eps}^2_{L_2(\Omega)}+\eps^2\norm{\omega_\eps}^2_{L_2(\Omega)}\right).\label{eqzetaeps}
\end{align}
Applying the techniques used in \cite{ELW12} to derive an estimate for the quantity corresponding to the right-hand side of \eqref{eqzetaeps}, we obtain from \eqref{eq5.10} and \eqref{eqzetaeps} that
$$\norm{\zeta_\eps}^2_{L_2(\Omega)}+\eps^2\norm{\omega_\eps}^2_{L_2(\Omega)}\leq C_8\eps^4.$$
This proves \eqref{eq5.6}. From \eqref{eq5.5} and \eqref{eq5.6} one concludes that $\norm{\partial_{z}\psi_\eps}_{W_2^1(\Omega)}\leq C_9\eps$ and \eqref{eq5.7} immediately follows from the embedding 7.56 in \cite[p.\ 217]{A75}. We set $\eps_*=\min\{\eps_1,\eps_2\}$ to complete our proof.
\endproof
Since we are interested in the limit $\eps\to 0$ of $(u_\eps,v_\eps,\phi_\eps)$, we have to guarantee that the maximal existence times $T_\eps>0$ do not converge to zero as $\eps\to 0$. Therefore, the following lemma generalizing \cite[Lem.~5.2]{ELW12} will be crucial.
\lem\label{lem_epstonull} There is $\tau=\tau(q,\lambda,\mu,\kappa)>0$ such that $\tau_\eps\geq\tau$ for all $\eps<\eps_*$. Moreover, there is $\Lambda=\Lambda(\kappa)>0$ such that $\tau_\eps=T_\eps=\infty$ for all $\eps<\eps_*$, provided $\lambda,\mu\in(0,\Lambda)$.
\endlem\rm
\proof Recalling the methods used to prove Proposition~\ref{propphig}, we see that, for fixed $2\sigma\in(\tfrac{1}{2}-\tfrac{1}{q},\tfrac{1}{2})$, there exists a positive constant $K_2(q,\kappa)$ such that
$$\norm{g_{\eps,i}(u_\eps(t),v_\eps(t))}_{W_{2,D}^{2\sigma}(I)}\leq K_2,\quad\forall t\in[0,\tau_\eps],\quad i=1,2.$$
With the aid of Duhamel's formula, see Section~\ref{sec_lwp}, we conclude that
\begin{align}
\norm{u_\eps(t)}_{W^2_{q,D}(I)}&\leq\frac{M}{\kappa}+\lambda MK_2\mathcal I(t),\nonumber\\
\norm{\hat v_\eps(t)}_{W^2_{q,D}(I)}&\leq\frac{M}{\kappa}+\mu MK_2\mathcal I(t),\nonumber\\
u_\eps(t)&\leq 0,\nonumber\\
v_\eps(t)&\geq -1,\nonumber\\
u_\eps(t)-v_\eps(t)&\geq 2\kappa-2(\lambda+\mu)MK_2\mathcal I(t).\nonumber
\end{align}
Let $m=\max\{\lambda,\mu\}$. As $\mathcal I(t)\to 0$ for $t\to 0$, there is $\tau=\tau(q,\lambda,\mu,\kappa)>0$ so that
$$\mathcal I(t)\leq\min\left\{\frac{1}{mK_2\kappa},\frac{(2M-1)\kappa}{4mM^2K_2}\right\},\quad\forall t\in[0,\tau].$$
It is clear that
$$\norm{u_\eps(t)}_{W^2_{q,D}(I)},\norm{\hat v_\eps(t)}_{W^2_{q,D}(I)}\leq\frac{1}{\kappa_1}$$
and
$$u_\eps(t)-v_\eps(t)\geq 2\kappa_1,\quad -1\leq v_\eps(t)<u_\eps(t)\leq 0,$$
for all $t\in[0,\tau]\cap[0,\tau_\eps]$. By the definition of $\tau_\eps$, we conclude $\tau_\eps\geq\tau$. Letting
$$\Lambda(\kappa):=\min\left\{\frac{1}{\kappa K_2\mathcal I(\infty)},\frac{(2M-1)\kappa}{4M^2K_2\mathcal I(\infty)}\right\}$$
and $\lambda,\mu\in(0,\Lambda(\kappa))$, we find that $T_\eps=\tau_\eps=\infty$, as was to be shown.
\endproof
We are now ready to present a proof of the following main theorem about convergence towards solutions of the small aspect ratio problem. Let $\mathbf 1_{A}$ denote the characteristic function of the set $A\subset\R^2$.
\thm Let $\lambda,\mu>0$, $q\in(2,\infty)$ and $\kappa\in(0,1/2)$ and let $(u_0,v_0)\in S_q(\kappa)$ satisfying $u_0\leq 0$ and $v_0\geq-1$ be given. For $\eps>0$, the unique solution to \eqref{originalproblem1}--\eqref{originalproblem8} with initial values $(u_0,v_0)$ obtained in Theorem~\ref{thm1} is denoted by $(u_\eps,v_\eps,\phi_\eps)$. The maximal interval of existence is $[0,T_\eps)$. Then there are $\tau>0$, $\eps_*>0$ and $\kappa_1\in(0,1/2)$ depending only on $q$ and $\kappa$ such that $T_\eps\geq\tau$ and $(u_\eps,v_\eps)(t)\in \overline S_q(\kappa_1)$ for all $(t,\eps)\in[0,\tau]\times(0,\eps_*)$. Moreover, the small aspect ratio model \eqref{sarm1}--\eqref{sarm7} has a unique solution $(u_*,v_*,\phi_*)$ so that
$$u_*,v_*\in C([0,\tau],W^2_{q}(I))\cap C^1([0,\tau],L_q(I)),$$
$-1\leq v_*(t)<u_*(t)\leq 0$ and $u_*(t)-v_*(t)\geq 2\kappa_1$ for $t\in[0,\tau]$, and such that, for a null sequence $(\eps_n)_{n\in\N}\subset(0,\eps_*)$,
\begin{align}
(u_{\eps_n},v_{\eps_n})&\to (u_*,v_*) &\text{in }C^{1-\theta}([0,\tau],W_q^{2\theta}(I)),\,\theta\in(0,1), \nonumber\\
\phi_{\eps_n}(t)\mathbf{1}_{\Omega_{u_{\eps_n}(t),v_{\eps_n}(t)}}&\to\phi_*(t)\mathbf{1}_{\Omega_{u_*(t),v_*(t)}} &\text{ in }L_2(I\times(-1,0)),\,t\in[0,\tau], \nonumber
\end{align}
as $n\to\infty$. Furthermore, there is $\Lambda(\kappa)>0$ such that, for $\lambda,\mu<\Lambda(\kappa)$, the statements of the theorem hold true for any $\tau>0$.
\endthm\rm
\proof Let $\tau$ and $\eps_*$ be as in Lemma~\ref{lem_epstonull}. Computing the $L_{q}(I)$-norm of \eqref{transfeq3} and \eqref{transfeq4} and using the reverse triangle inequality and the embedding $W^{1/2}_2(I)\hookrightarrow L_{2q}(I)$, we have, for any $t\in[0,\tau]$,
\begin{align}
&\norm{\partial_tu_\eps(t)}_{L_q(I)}-\norm{\partial_x^2u_\eps(t)}_{L_q(I)}\nonumber\\
&\qquad \leq K_3 \lambda\norm{\frac{1+\eps^2|\partial_xu_\eps(t)|^2}{(u_\eps(t)-v_\eps(t))^2}}_{L_\infty(I)}\norm{\partial_{z'}\tilde\phi_\eps(t,\cdot,1)}_{W^{1/2}_2(I)}^2\nonumber
\end{align}
and
\begin{align}
&\norm{\partial_tv_\eps(t)}_{L_q(I)}-\norm{\partial_x^2v_\eps(t)}_{L_q(I)}\nonumber\\
&\qquad \leq K_3 \mu\norm{\frac{1+\eps^2|\partial_xv_\eps(t)|^2}{(u_\eps(t)-v_\eps(t))^2}}_{L_\infty(I)}\norm{\partial_{z'}\tilde\phi_\eps(t,\cdot,0)}_{W^{1/2}_2(I)}^2,\nonumber
\end{align}
where $K_3>0$ is a constant. Since $\norm{\partial_{z'}\tilde\phi_\eps(t,\cdot,1)}_{W^{1/2}_2(I)}$ is bounded by a positive constant, cf.~Lemma~\ref{lemyyy}, we can proceed as in the proof of \cite[Thm.\ 1.4]{ELW12} to conclude from the boundedness of $u_\eps$ in $C([0,\tau],W^2_{q}(I))\cap C^1([0,\tau],L_q(I))$ that there exists a sequence $(u_{\eps_k})_{k\in\N}\subset\{u_\eps;\,\eps<\eps_*\}$, $\eps_k\to 0$, such that, for $k\to 0$,
$$
u_{\eps_k} \to u_* \text{  in  }\,C^{1-\theta}([0,\tau],W_q^{2\theta}(I))
$$
for some function $u_*\in C^{1-\theta}([0,\tau],W_q^{2\theta}(I))$ and $\theta\in\left(\frac{q+1}{2q},1\right)$. The boundedness of $\norm{\partial_{z'}\tilde\phi_{\eps_k}(t,\cdot,0)}_{W^{1/2}_2(I)}$, see again Lemma~\ref{lemyyy}, implies that $v_{\eps_k}$ is bounded in $C([0,\tau],W^2_{q}(I))\cap C^1([0,\tau],L_q(I))$ so that we may extract another null sequence $(\eps_{k_j})_{j\in\N}\subset(\eps_k)_{k\in\N}$ such that
$$
v_{\eps_{k_j}} \to v_* \text{  in  }\,C^{1-\theta}([0,\tau],W_q^{2\theta}(I)),
$$
for some function $v_*\in C^{1-\theta}([0,\tau],W_q^{2\theta}(I))$ and $\theta\in\left(\frac{q+1}{2q},1\right)$. As any subsequence of a convergent sequence is convergent with the same limit, we get that
\begin{align}
u_{\eps_{k_j}} \to u_* & \text{  in  }\,C^{1-\theta}([0,\tau],W_q^{2\theta}(I)) \text{ and}\nonumber\\
v_{\eps_{k_j}} \to v_* & \text{  in  }\,C^{1-\theta}([0,\tau],W_q^{2\theta}(I)) \nonumber
\end{align}
as $j\to\infty$. According to the continuous embedding $W_q^{2\theta}(I)\hookrightarrow W^1_\infty(I)$, we also have that
\begin{align}
u_{\eps_{k_j}} \to u_* & \text{  in  }\,C([0,\tau],W^1_\infty(I))\text{ and} \nonumber\\
v_{\eps_{k_j}} \to v_* & \text{  in  }\,C([0,\tau],W^1_\infty(I)). \nonumber
\end{align}
In view of the inequality \eqref{eq5.7} and the continuous embedding $W^{1/2}_2(I)\hookrightarrow L_{2q}(I)$ we observe that
$$\lim_{\eps\to 0}\sup_{t\in[0,\tau]}\norm{\left|\partial_{z'}\tilde\phi_\eps(t,\cdot,0)\right|^2-1}_{L_q(I)}=\lim_{\eps\to 0}\sup_{t\in[0,\tau]}\norm{\left|\partial_{z'}\tilde\phi_\eps(t,\cdot,1)\right|^2-1}_{L_q(I)}=0,$$
and conclude that
\begin{align}
&g_{\eps_{k_j}}(u_{\eps_{k_j}}(t),v_{\eps_{k_j}}(t))\to\frac{1}{(u_*-v_*)^2}(1,1)\,\text{  in  }\,C([0,\tau],L_q(I)), \nonumber\\
&u_*(t,x)\leq 0,\nonumber\\
&v_*(t,x)\geq-1,\nonumber\\
&u_*(t,x)-v_*(t,x)\geq 2\kappa_1,\nonumber\\
&(u_*,v_*)(t,\pm 1)=(0,-1)\text{ and}\nonumber\\
&(u_*,v_*)(0,x)=(u_0,v_0)(x)\text{ for }(t,x)\in[0,\tau]\times I.\nonumber
\end{align}
Moreover, for any $t\in[0,\tau]$, the left-hand side of the equation
$$u_{\eps_{k_j}}(t)=e^{-tA}u_0 -\lambda \int_0^t e^{-(t-s)A}g_{\eps_{k_j},1}(u_{\eps_{k_j}}(s),v_{\eps_{k_j}}(s))\,ds$$
converges to $u_*(t)$ while on the right-hand side the fact that $A$ generates the heat semigroup on $L_q(I)$ implies that
$$\int_0^t e^{-(t-s)A}g_{\eps_{k_j},1}(u_{\eps_{k_j}}(s),v_{\eps_{k_j}}(s))\,ds\to\int_0^t e^{-(t-s)A}\frac{1}{(u_*(s)-v_*(s))^{2}}\,ds,\quad j\to\infty.$$
Arguing similarly for $v_*$, we conclude that $u_*(t),\hat v_*(t)\in D(A)=W^2_{q,D}(I)$ and that $(u_*,v_*)(t)$ is the unique solution to \eqref{sarm2}--\eqref{sarm7} satisfying the properties stated in the theorem, for all $t\in[0,\tau]$.
%
%
%
Letting $$\phi_*=\frac{z-v_*}{u_*-v_*}$$ and using Lemma~\ref{lemyyy}, the proof is completed by similar arguments as in the proof of \cite[Thm.~2]{mk13}.
\endproof
\section{Discussion and Outlook}\label{sec_disc}
From the physical point of view, the effectiveness of a MEMS device is limited by the pull-in stability which corresponds to smash-up of both membranes in our idealized model. Intuitively, it is clear that this phenomenon occurs for large voltage values and thus, as the parameters $\lambda$ and $\mu$ are proportional to the square of the source voltage, cf.\ the Appendix, for large values of $\lambda$ and $\mu$. For the small aspect ratio model of a stationary MEMS device with a one-dimensional displacement of a single membrane suspended above a fixed ground plate, i.e.,
\beq\label{sarpullin}w_{xx}=\frac{\lambda}{(1+w)^2},\quad x\in(-1/2,1/2),\quad w(\pm 1/2)=0,\eeq
it is well-known that there is a threshold $\lambda_*$ such that for $0<\lambda<\lambda_*$ there exist two solutions $w_1(x;\lambda)$ and $w_2(x;\lambda)$ coalescing as $\lambda\to\lambda_*$ and there is no solution if $\lambda>\lambda_*$. Moreover, only one of the solutions in the small voltage regime is stable under perturbations; the other one is instable, cf.~\cite{BGP00,PT00}. Thus for this type of model, $\lambda_*$ corresponds indeed to the pull-in voltage. In \cite{BGP00}, the authors have computed the numerical value $\lambda_*=1.40001647737100$.

In Section~\ref{sec_lwp} we have first shown that, for any pair of sufficiently small parameters, there exists a solution $(u_\eps,v_\eps,\phi_\eps)$ to \eqref{originalproblem1}--\eqref{originalproblem8}, at least locally in time. Moreover, we have proven that there exists $m_1(\eps)>0$ such that $(u_\eps,v_\eps,\phi_\eps)$ is in fact a global solution, i.e., $T_\eps=\infty$, for $\lambda,\mu<m_1(\eps)$. In Section~\ref{sec_nonex} we have shown that there is $m_2(\eps)>0$ such that there is no global solution, i.e., $T_\eps<\infty$, for $\lambda>m_2(\eps)$ or $\mu>m_2(\eps)$. Note that our results do not provide information about the precise value of the pull-in voltage for this type of model. For instance, it is an open problem to find out whether the values $m_1$ and $m_2$ coincide or not. While Theorem~\ref{thm1} and Theorem~\ref{thmnonex} show that the sets
$$\mathcal S_\eps=\{(\lambda,\mu)\in(0,\infty)^2;\,\text{\eqref{originalproblem1}--\eqref{originalproblem8} has a global solution}\}$$
contain a neighborhood of zero in the relative topology of $(0,\infty)^2$ and are bounded by an $\eps$-dependent constant, we do not have further information on the structure of the $\mathcal S_\eps$; e.g., is not clear that $\mathcal S_\eps$ is the product of two intervals.

Moreover, it is not clear that $T_\eps<\infty$ implies that the membranes collide. One could also imagine that one component of the solution blows up in the corresponding $W_q^2$-norm; note that in the $W_q^2$-setting, $q\in(2,\infty)$, second order derivatives may become unbounded.

Concerning stability, we have already shown that there exists a steady state of \eqref{originalproblem1}--\eqref{originalproblem8} for sufficiently small parameters in \cite{mk13}. The present paper extends this result by proving uniqueness of the steady state (with first components in a set $S_q(\kappa)$) as well as its local asymptotic stability. It is an open problem whether there are other smooth branches of steady states emanating from $(\lambda,\mu)=(0,0)$ and what one can say about their stability or instability, cf.\ the discussion of the model \eqref{sarpullin} above.

Finally, the small aspect ratio limit has been discussed: We have first proven that the maximal existence times $T_\eps$ are bounded from below by a positive constant when sending $\eps\to 0$. Then refining the arguments of \cite[Thm.~1.4]{ELW12} and \cite[Thm.~2]{mk13} we have given a rigorous justification of the model \eqref{sarm1}--\eqref{sarm7} by proving convergence of the solutions $(u_\eps,v_\eps,\phi_\eps)$ towards a solution of \eqref{sarm1}--\eqref{sarm7} in the vanishing aspect ratio limit. Again, a cornerstone of our proof was to show that the arguments for the small aspect ratio limit of the stationary problem can be adopted for the evolution model and the $W_q^2$-setting.
\section{Appendix}
The mathematical model for an idealized electrostatic MEMS, considered in this paper, can be obtained as follows: There are two elastic membranes of length $\ell>0$ and width $w>0$ which are assumed to be perfect conductors and they should be fixed along their boundary so that their initial distance is $d>0$. We assume that a voltage $V_s$ is applied to the device so that an electric field with the potential $\psi$ sets up in the region bounded by the two membranes. Finally, let $\tilde u(\tilde x)$ and $\tilde v(\tilde x)$ denote the displacements of the membranes so that $(\tilde u,\tilde v)\equiv(0,-d)$ for $V_s=0$; see Fig.~\ref{figreal}.
We introduce coordinates $(\tilde x,\tilde y,\tilde z)\in\R^3$ so that the upper membrane is modeled by the set
$$\mathcal M_1=\{(\tilde x,\tilde y,\tilde z)\in\R^3;\;-\ell/2\leq\tilde x\leq \ell/2,\,-w/2\leq\tilde y\leq w/2,\,\tilde z=\tilde u(\tilde x)\}$$
and the second membrane corresponds to
$$\mathcal M_2=\{(\tilde x,\tilde y,\tilde z)\in\R^3;\;-\ell/2\leq\tilde x\leq \ell/2,\,-w/2\leq\tilde y\leq w/2,\,\tilde z=\tilde v(\tilde x)\}.$$
The region bounded by $\mathcal M_1$ and $\mathcal M_2$ is
$$\mathcal R=\{(\tilde x,\tilde y,\tilde z)\in\R^3;\;-\ell/2\leq\tilde x\leq \ell/2,\,-w/2\leq\tilde y\leq w/2,\,\tilde v(\tilde x)<\tilde z<\tilde u(\tilde x)\}.$$
The electrostatic potential satisfies the Laplace equation in $\mathcal R$ and we may choose $\mathcal M_2$ to be the set where $\psi=0$. Next, the potential on $\mathcal M_1$ is proportional to $V_s$ and there is a dimensionless function $f$ such that
\begin{align}
\frac{\partial^2\psi}{\partial\tilde x^2}+\frac{\partial^2\psi}{\partial\tilde y^2}+\frac{\partial^2\psi}{\partial\tilde z^2}=0 & \quad\text{in }\mathcal R,\label{Laporig}\\
\psi(\tilde x,\tilde y,\tilde z)=V_sf(\tilde u/d) & \quad\text{on }\mathcal M_1,\label{Lapbc1orig}\\
\psi(\tilde x,\tilde y,\tilde z)=0 & \quad\text{on }\mathcal M_2.\label{Lapbc2orig}
\end{align}
The function $f$ embodies the fact that the voltage drop across our device when embedded in a circuit may depend upon $\tilde u$, see also \cite{P01}. The ratio of the energy density of the electric field in $\mathcal R$ to the curvature of $\mathcal M_1$ and $\mathcal M_2$ is modeled by the surface tension coefficients $T_1,T_2>0$. With $\eps_0$ the permittivity of free space and $\eps_r$ the permittivity of the medium that fills $\mathcal R$ we thus have
\begin{align}
T_1\frac{\partial^2\tilde u}{\partial\tilde x^2} = \frac{1}{2}\eps_0\eps_r\left(\left(\frac{\partial\psi}{\partial\tilde x}\right)^2+\left(\frac{\partial\psi}{\partial\tilde y}\right)^2+\left(\frac{\partial\psi}{\partial\tilde z}\right)^2\right) & \quad\text{on }\mathcal M_1,\label{soduorig}\\
T_2\frac{\partial^2\tilde v}{\partial\tilde x^2} = -\frac{1}{2}\eps_0\eps_r\left(\left(\frac{\partial\psi}{\partial\tilde x}\right)^2+\left(\frac{\partial\psi}{\partial\tilde y}\right)^2+\left(\frac{\partial\psi}{\partial\tilde z}\right)^2\right) & \quad\text{on }\mathcal M_2.\label{sodvorig}
\end{align}
The sign in \eqref{sodvorig} is a consequence of the fact that both membranes should attract each other, and since $\mathcal M_1$ and $\mathcal M_2$ are fixed along the boundary, we have to impose the boundary conditions
\begin{align}
\tilde u(\ell/2)&=\tilde u(-\ell/2)=0,\label{bcuorig}\\
\tilde v(\ell/2)&=\tilde v(-\ell/2)=-d.\label{bcvorig}
\end{align}
The model \eqref{Laporig}--\eqref{bcvorig} is a free boundary problem, since the domain $\mathcal R$ and its boundary components $\mathcal M_1$ and $\mathcal M_2$ depend on the unknown functions $\tilde u$ and $\tilde v$ that also appear in the model equations. For $T_2\to\infty$, \autoref{sodvorig} takes the form $\tilde v''(\tilde x)=0$ and together with the boundary conditions \eqref{bcvorig} one immediately concludes that $\tilde v\equiv -d$. We thus recover the MEMS model with a fixed ground plate for $T_2\to\infty$ from the enhanced model presented here.\\
\begin{figure}[H]\label{figreal}
\begin{center}
\includegraphics[width=6.1cm]{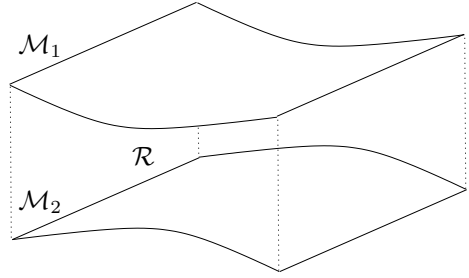}
\end{center}
\vspace{-.25cm}
\caption{Geometry of a three-dimensional MEMS device with two free boundaries $\mathcal M_1$ and $\mathcal M_2$ and a one-dimensional displacement.}
\label{figreal}
\begin{picture}(0,0)
\put(-83,95){$\mathcal M_1$}
\put(-83,36){$\mathcal M_2$}
\put(-40,53){$\mathcal R$}
\end{picture}
\end{figure}
To obtain dimensionless equations, we apply the transformation
$$\phi=\psi/V_s,\quad u=\tilde u/d,\quad v=\tilde v/d,\quad x=2\tilde x/\ell,\quad y=2\tilde y/w,\quad z=\tilde z/d$$
and introduce the dimensionless parameters
$$\eps=2d/\ell,\quad a=\ell/w,\quad \lambda=\frac{\eps_0\eps_rV_s^2\ell^2}{8T_1d^3},\quad\mu=\frac{\eps_0\eps_rV_s^2\ell^2}{8T_2d^3}$$
to rewrite the problem \eqref{Laporig}--\eqref{bcvorig} as
\begin{align}
0 & = \eps^2\frac{\partial^2\phi}{\partial x^2}+(\eps a)^2\frac{\partial^2\phi}{\partial y^2}+\frac{\partial^2\phi}{\partial z^2}, \\
\phi(x,y,u(x)) & = f(u),\\
\phi(x,y,v(x)) & = 0,\\
\frac{\partial^2u}{\partial x^2} & = \lambda \left.\left(\eps^2\left(\frac{\partial\phi}{\partial x}\right)^2+(\eps a)^2\left(\frac{\partial\phi}{\partial y}\right)^2+\left(\frac{\partial\phi}{\partial z}\right)^2\right)\right|_{z=u(x)},\label{bdry1}\\
\frac{\partial^2v}{\partial x^2} & = -\mu \left.\left(\eps^2\left(\frac{\partial\phi}{\partial x}\right)^2+(\eps a)^2\left(\frac{\partial\phi}{\partial y}\right)^2+\left(\frac{\partial\phi}{\partial z}\right)^2\right)\right|_{z=v(x)},\label{bdry2}\\
(u,v)(\pm 1) & = (0,-1).
\end{align}
The parameter $\eps$ is the small aspect ratio comparing gap size to device length, $a$ is the aspect ratio of the device itself and $\lambda$ and $\mu$ interrelate the strengths of the electrostatic and mechanical forces in the device. We did not discuss the case where the MEMS is embedded into a circuit so that we can set $f\equiv 1$ in the following. Moreover, the fact that we have assumed a one-dimensional displacement of the membranes motivates to assume that $\phi$ is a function of $x$ and $z$ only.

Assuming that $u$ and $v$ are functions of time $\tilde t$ and applying Newton's second law on both membranes, we obtain that the sum of all forces equals $\rho_1\delta_1\frac{\partial^2u}{\partial{\tilde t}^2}$ and $\rho_2\delta_2\frac{\partial^2v}{\partial{\tilde t}^2}$, where $\rho_1,\rho_2$ and $\delta_1,\delta_2$ denote the mass density per unit volume of the membranes and the membrane thicknesses. With a damping force term of the form $-\sigma\frac{\partial u}{\partial{\tilde t}}$ and $-\sigma\frac{\partial v}{\partial{\tilde t}}$ respectively, \autoref{bdry1} and \autoref{bdry2} take the form
\begin{align}
\rho_1\delta_1\frac{\partial^2 u}{\partial{\tilde t}^2} + \sigma\frac{\partial u}{\partial{\tilde t}} & = \frac{\partial^2u}{\partial x^2}-\lambda\left.\left(\eps^2\left(\frac{\partial\phi}{\partial x}\right)^2+\left(\frac{\partial\phi}{\partial z}\right)^2\right)\right|_{z=u(x)}, \\
\rho_2\delta_2\frac{\partial^2 v}{\partial{\tilde t}^2} + \sigma\frac{\partial v}{\partial{\tilde t}} & = \frac{\partial^2 v}{\partial x^2}+\mu\left.\left(\eps^2\left(\frac{\partial\phi}{\partial x}\right)^2+\left(\frac{\partial\phi}{\partial z}\right)^2\right)\right|_{z=v(x)} .
\end{align}
Setting $t=\tilde t/\sigma$, $\gamma_1=\sqrt{\rho_1\delta_1}/\sigma$ and $\gamma_2=\sqrt{\rho_2\delta_2}/\sigma$, one finally has
\begin{align}
\gamma_1^2\frac{\partial^2 u}{\partial{t}^2} + \frac{\partial u}{\partial{t}} & = \frac{\partial^2u}{\partial x^2}-\lambda\left.\left(\eps^2\left(\frac{\partial\phi}{\partial x}\right)^2+\left(\frac{\partial\phi}{\partial z}\right)^2\right)\right|_{z=u(x)}, \\
\gamma_2^2\frac{\partial^2 v}{\partial{t}^2} + \frac{\partial v}{\partial{t}} & = \frac{\partial^2 v}{\partial x^2}+\mu\left.\left(\eps^2\left(\frac{\partial\phi}{\partial x}\right)^2+\left(\frac{\partial\phi}{\partial z}\right)^2\right)\right|_{z=v(x)} .
\end{align}
In this paper, we have assumed that $\gamma_1,\gamma_2\ll 1$ to obtain the problem \eqref{originalproblem1}--\eqref{originalproblem8} with parabolic equations on the free boundaries. To study the hyperbolic-elliptic free boundary problem with $\gamma_1,\gamma_2>0$ is a task for further research.

\end{document}